\documentclass[psamsfonts]{amsart}
\usepackage{amsmath,amsthm,amssymb,amsxtra,latexsym}

\usepackage{pifont}
\usepackage{enumerate}

\usepackage{xypic}

\def\into{{\hookrightarrow}}
\def\xto{\xrightarrow}

\newcommand{\fS}{{\mathfrak S}}
\newcommand{\fX}{{\mathfrak X}}

\newcommand{\calf}{{\mathcal F}}

\newcommand{\calh}{{\mathcal H}}
\newcommand{\cali}{{\mathcal I}}
\newcommand{\calj}{{\mathcal J}}

\newcommand{\call}{{\mathcal L}}
\newcommand{\calm}{{\mathcal M}}
\newcommand{\caln}{{\mathcal N}}
\newcommand{\calo}{{\mathcal O}}

\newcommand{\calt}{{\mathcal T}}

\newcommand{\calos}{{\calo_S}}

\newcommand{\calox}{{\calo_X}}

\newcommand{\caloz}{{\calo_Z}}

\def\cHom{{\mathcal H}om}

\def\cDer{{\mathcal D}er}

\newcommand{\bbbb}{{\mathbb B}}
\newcommand{\bbbc}{{\mathbb C}}

\newcommand{\bbbl}{{\mathbb L}}

\newcommand{\bbbm}{{\mathbb M}}
\newcommand{\bbbp}{{\mathbb P}}

\newcommand{\bbbv}{{\mathbb V}}
\newcommand{\bbbz}{{\mathbb Z}}

\newcommand{\vp}{\varphi}

\DeclareMathOperator{\codim}{codim}

\DeclareMathOperator{\coker}{coker}

\DeclareMathOperator{\Der}{Der}
\DeclareMathOperator{\depth}{depth}

\DeclareMathOperator{\End}{ End}

\DeclareMathOperator{\Ext}{Ext}

\DeclareMathOperator{\grad}{grad}

\DeclareMathOperator{\Hilb}{Hilb}
\DeclareMathOperator{\Hom}{Hom}

\DeclareMathOperator{\Image}{Im}
\DeclareMathOperator{\jac}{jac}

\DeclareMathOperator{\Sing}{Sing}

\newcommand{\CC}{{\mathbb C}}
\newcommand{\HH}{{\mathbb H}}
\newcommand{\PP}{{\mathbb P}}
\newcommand{\QQ}{{\mathbb Q}}

\newcommand{\ZZ}{{\mathbb Z}}

\newcommand{\y}{{s}}

\newcommand{\GG}{{\mathbb G}}
\newcommand{\GGdual}{\check{\GG}}
\newcommand{\GGdualblowup}{\tilde{\check{\GG}}}
\newcommand{\PPdual}{\check{\PP}}
\newcommand{\PPdualblowup}{\tilde{\check{\PP}}}
\newcommand{\Wperp}{W_\perp}
\newcommand{\Uperp}{U_\perp}
\newcommand{\Xdual}{\check{X}}
\newcommand{\Yperp}{Y_\perp}
\newcommand{\yperp}{y_\perp}
\newcommand{\PPperp}{\PP_\perp}
\newcommand{\Edual}{\check{E}}
\newcommand{\Eperp}{E_\perp}
\newcommand{\Bperp}{B_\perp}

\newcommand{\sF}{{\mathcal F}}
\newcommand{\sO}{{\mathcal O}}
\newcommand{\sQ}{{\mathcal Q}}
\newcommand{\sL}{{\mathcal L}}
\newcommand{\Jac}{{\rm Jac}}

\theoremstyle{definition}
\newtheorem{defn}{Definition}[section]
\newtheorem{definition}[defn]{Definition}

\newtheorem{example}[defn]{Example}

\newtheorem{sit}[defn]{}
\newtheorem{remark}[defn]{Remark}

\theoremstyle{plain}
\newtheorem{prop}[defn]{Proposition}
\newtheorem{proposition}[defn]{Proposition}
\newtheorem{theorem}[defn]{Theorem}

\newtheorem{lemma}[defn]{Lemma}
\newtheorem{cor}[defn]{Corollary}

\theoremstyle{remark}

\begin{document}
\title[Free Divisors] {Low-dimensional Singularities with Free Divisors as Discriminants}

\author[R.-O.~Buchweitz]{Ragnar-Olaf Buchweitz}
\address{Dept.\ of Computer and Mathematical Sciences, University of Tor\-onto at Scarborough, Tor\-onto, Ont.\ M1A 1C4, Canada}
\email{ragnar@math.utoronto.ca}

\author[W.~Ebeling]{Wolfgang Ebeling}
\address{Institut f\"ur Algebraische Geometrie, Leibniz Universit\"at Hannover, Postfach 
6009, D-30060 Hannover, Germany}
\email{ebeling@math.uni-hannover.de}

\author[H.-C.~Graf~von~Bothmer]{Hans-Christian Graf von Bothmer}
\address{nstitut f\"ur Algebraische Geometrie, Leibniz Universit\"at Hannover, Postfach 
6009, D-30060 Hannover, Germany}
\email{bothmer@math.uni-hannover.de}

\thanks{The authors were partly supported by the DFG Schwerpunkt ``Global Methods in Complex Geometry". The first author was also partly supported by NSERC grant
3-642-114-80 and wishes to thank his alma mater, the University of Hannover, as well as S.-O.~Buchweitz-Klings\"ohr and D.~Klings\"ohr for their hospitality during the preparation of this work.}


\begin{abstract}
We present versal complex analytic families, over a smooth base
and of fibre dimension zero, one, or two, where the discriminant constitutes
a free divisor. These families include finite flat maps, versal deformations of reduced curve singularities, and versal deformations of Gorenstein surface singularities in $\CC^5$. It is shown that such free divisors often admit a ``fast normalization'', obtained by a single application of the Grauert-Remmert normalization algorithm. For a particular Gorenstein surface singularity in $\CC^5$, namely the simple 
elliptic singularity of type $\widetilde A_{4}$, we exhibit an explicit 
discriminant matrix and show that the slice of the discriminant for a fixed $j$--invariant is the 
cone over the dual variety of an elliptic curve.
\end{abstract}

\maketitle

{\footnotesize\tableofcontents}
\section*{Introduction}
One of the remarkable results in complex singularity theory is that the 
discriminant in the versal deformation of any isolated complete 
intersection singularity is a {\em free divisor\/}, a highly singular 
hypersurface in the ambient smooth base that admits though a 
compact representation as determinant of any {\em discriminant matrix\/} 
expressing a basis of the liftable vector fields in terms of the coordinate 
vector fields; see \cite{Lo4} or Section 2 below.

Variants on this result show the freeness of the discriminant  
in the base of a versal deformation in a number of further cases,
for example space-curve singularities (van Straten \cite{vSt}), functions on space curves 
(Goryunov  \cite{Gor}, Mond-van Straten \cite{MvSt}), or representation varieties of quivers 
(Buchweitz-Mond \cite{BMo}). 
J.~Damon gives in \cite{Da2} a broad survey of why and how free divisors appear 
frequently in the theory of discriminants and bifurcations.

Here we present further versal complex analytic families, over a smooth base
and of fibre dimension zero, one, or two, where the discriminant constitutes
yet again a free divisor. 

While we consider the question in more generality to deduce sufficient criteria, 
showing along the way why such free divisors often admit a ``fast normalization'', 
obtained by a single application of the Grauert-Remmert normalization algorithm,
the new explicit instances found pertain to 
\begin{itemize}
\item finite flat maps, thus the case of relative dimension zero,
where we characterize freeness of the discriminant completely and 
show, for example, that the discriminant in the Hilbert scheme or Douady space
of points on a smooth complex surface is a free divisor,
\item reduced curve singularities, where we recover not only the result on 
space-curve singularities  due to D.~van Straten \cite{vSt}, but extend it to
include all reduced, smoothable, and unobstructed Gorenstein 
curve singularities, and
\item smoothable Gorenstein surface singularities with reflexive conormal module,
thus including all Gorenstein surface singularities in $\bbbc^{5}$.
\end{itemize}

This list includes in particular versal deformations of any isolated 
Gorenstein singularity of dimension at most two that is linked to a 
complete intersection. 
More generally, we start with an isolated Cohen-Macaulay singularity, 
of arbitrary dimension, whose conormal module becomes Cohen-Macaulay 
when twisted with the canonical module.
J.~Herzog \cite{Her} showed that such a singularity is unobstructed, 
whence any versal deformation has a smooth base, and that the module 
of vector fields on the total space of a versal deformation is Cohen-Macaulay.
Assuming that the singular locus of the total space is not dense in the critical locus
of the versal morphism, the discriminant is necessarily a divisor, and freeness becomes equivalent to sufficient depth of the module of {\em vertical vector fields\/} in the versal deformation. That latter property turns out to be more and more elusive as the 
dimension of the singularity increases, and it is not quite clear what the 
actual reach of the criterion is. 
However, in small enough dimension, the necessary depth can be verified easily,
as we show.

Once  these rather abstract results of homological nature are presented, 
we turn to the concrete description of some of the free divisors that arise 
and point out how they are related to dual varieties of smooth 
projective varieties. The avatar here is the classical result that the dual 
variety to the rational normal curve yields the discriminant in the space 
of univariate polynomials. In particular, building on work by many authors, 
mainly H.~Pinkham \cite{Pi1, Pi2}, E.~Looijenga \cite{Lo, Lo2, Lo3}, K.~Saito \cite{Sa2}, 
K.~Wirthm\"uller \cite{Wi}, J.-Y.~M\'erindol \cite{M}, B.~Dubrovin \cite{Du}, and M.~Bertola \cite{Be2, Be3, Be4}, we put a capstone on 
the structure of the discriminant in the semi-universal deformation of a simple 
elliptic surface singularity of type $\widetilde A_{4}$, exhibiting an explicit 
discriminant matrix and showing that its slice for a fixed $j$--invariant is the 
cone over the dual variety of the dual elliptic curve. Note that according to 
\cite{Sa1, Pi1, Hu} a surface singularity of type $\widetilde A_{4}$ is the only 
instance of a simple elliptic singularity, where the deformation theory is 
unobstructed but the singularity itself is not a complete intersection, as it is rather given minimally by the Pfaffians of an alternating $(5\times 5)$--matrix, thus Gorenstein of codimension $3$ and therefore (only) linked to a complete intersection.

Specifically, the paper is organized as follows. In Section 1, we recall basic properties
of free divisors and indicate how the appearance of free divisors in versal deformations
can be explained in terms of the Kodaira-Spencer sequence.
Section 2 contains the main result, a general criterion for the discriminant
in the base of a versal deformation to be a free divisor. Section 3 investigates the case
of finite flat maps with applications to the Hilbert scheme of points, Section 4 applies
the main result to families of curves and surfaces. Section 5 contains a review of the
classical discriminant of a polynomial and shows that Arnol'd's description is equivalent
to the Bezout representation of the dual variety of the rational normal curve.
In Section 6 we show that the analogue of the Bezout formula found by M.~Bertola can be used to give a determinantal expression for the dual variety of the elliptic normal curve. 
After reviewing basic material on the deformation theory of simple elliptic 
singularities of type $\widetilde A_{4}$, we present in Section~Ê\ref{section:A4}
an explicit discriminant matrix and relate it to the dual variety of an elliptic curve.

{\em Acknowledgements:\/} We wish to thank K.~Wirthm\"uller and the
late P.~Slodowy for valuable discussions and pointers to the literature. For the preparation of Section~\ref{section:A4} the computer algebra systems {\sc Macaulay2} \cite{M2} and {\sc Maple}${}^{\rm TM}$ were used. 

\section{Background on Deformations and Free Divisors } 
We review here the definition and basic properties of free divisors. Further details
and a broad survey of possible generalizations are given in \cite{Da1}.

\begin{sit}
Let $S$ be a connected complex manifold and $D\subset S$ a reduced hypersurface 
with defining ideal sheaf $\cali_{D}\subseteq \calos$.
The sheaf $\cDer_{S}(-\log D)$ of {\em logarithmic vector fields\/}   along $D$ has as its local sections those vector fields $\chi \in \cDer_{S}$ such that $\chi(\cali_{D})\subseteq \cali_{D}$, or, equivalently, such that $\chi$ is tangent to $D$ at its regular points. It is clearly
an $\calos$-submodule and a sheaf of complex Lie subalgebras of $\cDer_{S}$.
\end{sit}

\begin{definition}
The reduced hypersurface $D\subset S$ is a {\em free divisor\/} at $s\in S$, if $\cDer_{S}(-\log D)$ is a locally free $\calos$-module at $s$, necessarily then of the same rank as $\cDer_{S}$, equal to the dimension of $S$.
\end{definition}

The concept of free divisors was identified by K.\,Saito in \cite{Sai} and he stated there
the following criterion, nowadays usually named after him:

\begin{proposition}{\em ({\sc Saito's Criterion})}
\label{scrit} 
The hypersurface $D\subset S$ is a free divisor at a point $s\in S$ if and 
only if there are germs  of vector fields $\chi_1,\ldots,\chi_n \in\cDer_{S}(-\log D)_{s}$,
such that the determinant of the matrix of coefficients of these germs
with respect to some, or any, 
$\calo_{S,s}$-basis of $\cDer_{S,s}$, is a reduced equation for $D$ at $s$. 
In this case, $\chi_1,\ldots, \chi_n$ form a basis for the $\calo_{S,s}$--module
$\cDer_{S}(-\log D)_{s}$.\qed
\end{proposition}

Any {\em discriminant matrix\/} that describes the inclusion $\cDer_{S}(-\log \Delta)\subseteq \cDer_{S}$, yields thus through its determinant a compact presentation of a defining equation of $D$. 

\begin{sit}
\label{sit:Aleks}
Free divisors are very special hypersurfaces. If not smooth, they 
are ``{\em maximally singular\/}'' in that the singular locus is a 
Cohen-Macaulay subspace of codimension one in $D$.
Equivalently, the Jacobi ideal $\jac_{D,s}\subseteq \calo_{D,s}$ 
of the free divisor, naturally isomorphic to the cokernel 
$\cDer_{S}/ \cDer_{S}(-\log D)$, 
is a {\em maximal Cohen-Macaulay module\/} of rank one on $D$, and any 
presentation matrix of it yields a discriminant matrix. 

Conversely,  as observed by A.~G.~Alexandrov \cite{Al1,Al2}; see also \cite{Da1}; 
a reduced complex hypersurface whose Jacobi ideal constitutes a maximal Cohen-Macaulay 
module on it is necessarily a free divisor.
\end{sit}

These characterizations can be formulated in concrete algebraic terms, without 
explicit reference to vector fields, just in terms of the Taylor series of a locally defining equation.

\begin{proposition}
A (formal) power series $f\in P:=\bbbc[\![z_{1},...,z_{n}]\!]$ defines a (formal) 
free divisor if it is reduced, that is, squarefree, and there is an 
$(n{\times}n)$--matrix $A$ with entries from $P$ such that
$$
\det A = f\quad\text{and}\quad 
(\nabla f)A\equiv (0,...,0) \bmod f\,,
$$
where $\nabla f =
\left(
\frac{\partial f}{\partial z_{1}},\ldots,\frac{\partial f}{\partial z_{n}}
\right)$ is the gradient of $f$, and the last condition just expresses that 
each entry of the (row) vector $(\nabla f)A$ is divisible by $f$ in $P$. 
The columns of $A$ can then be viewed as the coefficients of a basis, with 
respect to the partial derivatives $\partial/\partial z_{i}$,  of the 
(formal) logarithmic vector fields along the divisor $f=0$, and the cokernel 
of $A$ is naturally isomorphic to the Jacobi ideal of $f$ in $P/(f)$.\qed
\end{proposition}

\begin{sit}
A normal crossing divisor is a free divisor in any dimension, 
and in the theory of hyperplane arrangements many {\em free arrangements\/} 
have been constructed by combinatorial means; see e.g. \cite[Ch. 4]{OT}.
We note in passing that it is one of the major outstanding problems in that 
theory whether freeness is solely a combinatorial property. 
Going beyond such unions of locally linear spaces, while any 
reduced plane curve is a free divisor in its ambient plane, the concept
becomes quite elusive in higher dimensions. Even free surfaces, in
an ambient smooth space of dimension $3$, are not classified yet, 
only a zoo of rather few specimen is known; see \cite{Da2} for some of those. 
\end{sit}

We now turn to the appearence of free divisors as discriminants in 
versal deformations. To this end, we next review some standard notation 
and results pertaining to (versal) deformations of  complex 
spaces or germs thereof.

\begin{sit} Recall the notion of {\em tangent cohomology\/}, 
the groups that control infinitesimal deformations. For
any morphism $f:X\to S$ of complex analytic germs or spaces, let
$\bbbl_{X/S}\in D^{-}(X)$ denote an (analytic) {\em  cotangent 
complex\/} of $X$ over
$S$, or rather of $f$. Up to isomorphism, such cotangent complex is a 
well defined object in the indicated derived category,
see \cite{Fle} or \cite{BFl1, BFl}.

The cohomology group 
$T^{i}_{X/S}(\calm):=H^{i}(\Hom_{\calox}(\bbbl_{X/S},\calm))$  is the 
{\em $i^{th}$ tangent cohomology\/}
of the morphism $f$ with values in the $\calox$--module $\calm$. We 
abbreviate as usual
$T^{i}_{X/S}:= T^{i}_{X/S}(\calox)$. Similarly, one defines the 
{\em tangent cohomology sheaves\/}
$\calt^{i}_{X/S}(\calm):=\calh^{i}(\cHom_{\calox}(\bbbl_{X/S},\calm))$. 
If $f$ is a morphism of analytic germs and $\calm$ is
$\calox$--coherent, then each
$\calt^{i}_{X/S}(\calm)$ is a coherent $\calox$--module as well.
Moreover, tangent cohomology localizes in that 
$\calt^{i}_{X/S}(\calm)_{x}\cong T^{i}_{(X,x)/(S,f(x))}(\calm_{x})$
as $\calo_{X,x}$--modules for any point $x$ on $X$. Therefore,
we may, and will, use the sheaf or module-theoretic language 
interchangeably throughout for morphisms of germs.
\end{sit}

\begin{sit}
Note, in particular, that for $\calm$ coherent,
$T^{0}_{X/S}(\calm)\cong \Der_{X/S}(\calm)$,  the 
$\calox$--module of $\calos$--linear derivations on $\calox$  with
values in $\calm$, and that $T^{1}_{X/S}(f^{*}\caln)$, for a flat 
morphism $f$ and  a coherent $\calos$--module $\caln$,
parametrizes the isomorphism classes of extensions  of $X$ by 
$f^{*}\caln$ over the trivial extension of $S$ by $\caln$.
\end{sit}

\begin{sit} A morphism $f:X\to S$ is {\em versal\/} at $s\in S$, if 
it is flat and induces a formally versal deformation of
the fibre $X(s)=f^{-1}(s)$ at $s$ in $S$. It is {\em versal\/} if it 
is so at every point in $S$. Versality is an open
property on $S$ by \cite{Fle}; see also \cite{BFl}.

If, for a complex germ $X_{0}$, the first tangent cohomology $T^{1}_{X_{0}}$ is a 
finite dimensional vector space, then $X_{0}$ admits a versal deformation.  
If $T^{2}_{X_{0}}=0$, then the deformation theory of $X$ is unobstructed, and 
so the base of any versal deformation, if it exists, is necessarily smooth.
\end{sit}

\begin{defn}
\label{item:liftable}  Let $f:(X,0)\to (S,0)$ be a flat morphism of 
analytic germs\footnote{Most  base points of germs are
called ``$0$''. They are usually suppressed from the  notation.}. For 
any coherent $\calos$--module $\caln$, set
\begin{align*} T^{0}_{X\to S}(\caln)\cong &\{(D,D')| D\in 
\Der(\calos,\caln),D'\in\Der(\calox,
  f^{*}\caln),\\&\qquad D'\circ f = D\otimes_{\calos}1:\calos\to
  \caln\otimes_{\calos}\calox = f^{*}\caln\}\,.
\end{align*} In other words, this $\calos$--module consists of {\em 
compatible pairs of vector fields\/}, one on $S$ with
values in $\caln$, the other on $X$ with values in $f^{*}\caln$.

For $\caln=\calos$, the $\calos$--module $T^{0}_{X\to S}:=T^{0}_{X\to 
S}(\calos)$ carries naturally a structure of complex Lie
algebra with respect to the componentwise bracket of vector fields.

The projection $p_{1}$ from $T^{0}_{X\to S}$ to $T^{0}_{S}\cong 
\Der_{S}$ is a homomorphism of such Lie algebras. Its
image $\call\subseteq T^{0}_{S}$ is a complex Lie subalgebra and 
$\calos$--submodule, called the submodule of {\em liftable\/}
vector fields on $S$.

Similarly, the image of $T^{0}_{X\to S}(\calos)$ under the projection 
$p_{2}$  to $T^{0}_{X}(f^{*}\calos)\cong \Der_{X}$
consists of what Arnol'd calls the  {\em lowerable\/} vector fields on 
$X$. Those form again a complex Lie subalgebra,  but
only an $\calos$--submodule of $T^{0}_{X}$.
\end{defn}

\begin{sit}
\label{sit:basic} The relevance of the $\calos$--module $T^{0}_{X\to 
S}(\caln)$ of  compatible vector fields is its relation
to the {\em Kodaira-Spencer map\/}, a relation described through 
the following commutative 
diagram, with right column the exact {\em Zariski-Jacobi sequence\/} 
of the tangent cohomology $\calox$--modules associated to $f$ and $f^{*}\caln$,
and left column an exact sequence of $\calos$--modules\footnote{
To exhibit indeed all modules in the  left 
column as $\calos$--modules, one should write 
$f_{*}T^{1}_{X/S}(f^{*}\caln)$ there instead of $T^{1}_{X/S}(f^{*}\caln)$.
However, to keep notation at bay, we allow for this abuse of notation.},
that essentially reflects the description of $T^{0}_{X\to S}(\caln)$ as a 
fibred product. The horizontal arrows represent maps that are linear 
over $f$.
  \begin{align}
\label{diag:ZJKSseq}
\xymatrix{ 0\ar[d]&0\ar[d]\\ 
T^{0}_{X/S}(f^{*}\caln)\ar@{=}[r]\ar[d]&T^{0}_{X/S}(f^{*}\caln)\ar[d]\\ 
T^{0}_{X\to
S}(\caln)\ar[r]^{p_{2}}\ar[d]_{p_{1}}&T^{0}_{X}(f^{*}\caln)\ar[d]\\
T^{0}_{S}(\caln)\ar[r]^{can}\ar[d]_{\delta_{X/S}^{\caln}}&T^{0}_{S}(f_{*}f^{*}\caln)\ar[d]\\
T^{1}_{X/S}(f^{*}\caln)\ar@{=}[r]&T^{1}_{X/S}(f^{*}\caln)\ar[d]\\
&T^{1}_{X}(f^{*}\caln)\ar[d]\\
&T^{1}_{S}(f_{*}f^{*}\caln)\ar[d]\\
&\vdots }
\end{align}
The map $\delta_{X/S}^{\caln}$ is the {\em Kodaira-Spencer map\/} 
associated to $f$ and $\caln$, while the map
labeled $can$(-onical) is induced from the natural $\calos$--homomorphism
$\caln\to f_{*}f^{*}\caln\cong \caln\otimes_{\calos}\calox, n\mapsto 
n\otimes 1$,  the unit of the adjunction.
\end{sit}

\begin{sit}
\label{sit:facts} Note the following facts:
\begin{enumerate}
\item
\label{versalitycrit}  The {\em versality criterion}, see \cite{Fle} 
or \cite{BFl}, states that for any coherent
$\calos$--module $\caln$ the support of the cokernel of the 
Kodaira-Spencer  map $\delta_{X/S}^{\caln}$ is contained in the
locus of $S$ where $f$ is not {\em versal}. In particular, if $f$ is 
versal, then $\delta_{X/S}^{\caln}$ is surjective for any
coherent $\calos$--module $\caln$.

\item
\label{shortKS}
In particular, for $f$ versal and $\caln=\calos$, the left column in
(\ref{diag:ZJKSseq}) above yields the following short exact sequence of
$\calos$--modules that exhibits the liftable vector fields $\call$ from
\ref{item:liftable} as kernel of the Kodaira-Spencer map,
\begin{align*}
0\to \call \xto{\quad} \Der_{S}\xto{\delta_{X/S}} T^{1}_{X/S}\to 0\,.
\end{align*}

\item If $S$ is {\em smooth}, the functors $T^{i}_{S}(?)$ vanish on any
$\calos$--module for $i\neq 0$.
\item
\label{rigid}  Note as well that $T^{1}_{X}(f^{*}\calos)\cong 
T^{1}_{X}$ vanishes if,  and only if, $X$ is {\em rigid\/}.
\end{enumerate}
\end{sit} These facts have the following immediate consequences.
\begin{cor} If $f:X\to S$ is versal, then the $\calos$--module 
$f_{*}T^{1}_{X/S}(f^{*}\caln)$ is {\em coherent\/} along with
$\caln$.
\end{cor}
\begin{proof} Indeed, as $\calos$--module, $T^{1}_{X/S}(f^{*}\caln)$ 
is the homomorphic image of $T^{0}_{S}(\caln)$ under the
Kodaira-Spencer map by the versality criterion 
\ref{sit:facts}(\ref{versalitycrit}).
\end{proof}

\begin{cor}
\label{cor:rigid} If $f:X\to S$ is versal with $S$ smooth, then 
$T^{1}_{X}(f^{*}\caln)=0$ for any coherent $\calos$--module
$\caln$. In particular, $X$ is rigid.
\end{cor}
\begin{proof} As $f$ is versal, in the diagram (\ref{diag:ZJKSseq}) 
the Kodaira-Spencer map
$\delta^{\caln}_{X/S}$ is surjective by \ref{sit:facts}( 
\ref{versalitycrit}),  whence the map $T^{0}_{S}(f_{*}f^{*}\caln)\to
T^{1}_{X/S}(f^{*}\caln)$ must be surjective too.  As 
$T^{1}_{S}(f_{*}f^{*}\caln)=0$ by smoothness of $S$, it follows that
$T^{1}_{X}(f^{*}\caln)=0$ as claimed. The final assertion follows 
then  from \ref{sit:facts}(\ref{rigid}).
\end{proof}

\section{Free Divisors in Versal Deformations}
We are now mainly interested in the case where $f$ is a versal morphism of germs with
$S$ smooth and  where the  $\calos$--module $\call$ of liftable vector 
fields as in \ref{item:liftable} is {\em free\/}. From \ref{sit:facts}(\ref{shortKS}) one
obtains immediately the following equivalent characterizations,  and much 
of the subsequent work will be to establish and verify
other, more manageable characterizations, such as the one exhibited 
in \ref{lemma:depth} below.
\begin{lemma}
\label{lem:depthT1} If $f:X\to S$ is versal with $S$ smooth, then the 
following conditions are equivalent.
\begin{enumerate}
\item The submodule $\call\subseteq T^{0}_{S}$ of liftable vector 
fields is a  {\em free\/} $\calos$--module.
\item The $\calos$--module $f_{*}T^{1}_{X/S}$ is of projective dimension 
at most $1$.
\item The $\calox$--module $T^{1}_{X/S}$ is of depth at least $\dim S - 1$.
\end{enumerate}
\qed
\end{lemma}

\begin{cor}
\label{cor:gensmooth} If the equivalent conditions of 
\ref{lem:depthT1} are satisfied and if furthermore
$f|_{\Sing X}:\Sing X\to S$ is {\em not dominant\/}, then
\begin{enumerate}
\item The free $\calos$--module $\call$ of liftable vector fields is
of rank equal to $\dim S$.
\item The {\em zeroth Fitting ideal\/} $\calf_{0}f_{*}T^{1}_{X/S}\subseteq \calos$ 
of $T^{1}_{X/S}$ as $\calos$--module, is
{\em principal\/}, generated by the determinant $\Delta$ of any 
matrix representing the inclusion $\call\subseteq T^{0}_{S}$
of free $\calos$--modules of the same rank.
\item The $\calos$--module $f_{*}T^{1}_{X/S}$ is {\em maximal 
Cohen-Macaulay\/} on the hypersurface, or divisor,
$V(\calf_{0}f_{*}T^{1}_{X/S})\subseteq S$.
\item
\label{item:empty}  The support of $T^{1}_{X/S}$ is {\em empty\/}, 
equivalently,  $T^{1}_{X/S}=0$ or, also,
$\calf_{0}f_{*}T^{1}_{X/S}= \calos$, if, and only if, $f$ is {\em smooth\/}.
\end{enumerate}
\end{cor}

\begin{proof} The theorem on generic smoothness implies that over 
some Zariski-open  dense subset $S \setminus f(\Sing X)$, 
the morphism $f$ is smooth,  whence $f_{*}T^{1}_{X/S}$ is supported in 
$S\setminus U$. This implies that
$f_{*}T^{1}_{X/S}$ is a torsion $\calos$--module and forces the free 
$\calos$--module
$\call$ to be of the same rank as $T^{0}_{S}=\Der(\calos)$, which 
equals $\dim S$.  That the support of $f_{*}T^{1}_{X/S}$ is then
given by the indicated  Fitting ideal is nothing but Cramer's rule. 
The assertion that $f_{*}T^{1}_{X/S}$ is maximal Cohen-Macaulay
on its support in $S$ just restates that $T^{1}_{X/S}$, if not zero, 
is of projective dimension $1$ over $\calos$.

The module $T^{1}_{X/S}$ is zero if, and only if, its support is 
empty,  which happens thus if, and only if, every vector
field on $S$ can be lifted,  and that property in turn is equivalent 
to $f$ being smooth.
\end{proof}

Note that $\calf_{0}T^{1}_{X/S} \subseteq \calox$ defines the natural analytic structure on the critical locus $C(f)$ of $f$ in $X$, see e.g.\ \cite[4.A]{Lo4}.

\begin{defn} With assumptions and notation as in \ref{cor:gensmooth}, 
a matrix representing  the containement $\call\subseteq
T^{0}_{S}$ is called a {\em discriminant matrix\/}  for $f$ and its 
determinant $\Delta$ constitutes a {\em canonical
equation\/} for the resulting  divisor. The zero-set $V(\Delta)$ will 
also be called the (homological)  {\em  discriminant\/}
of $f$. It is not necessarily reduced.

Note that the discriminant might be empty.  By \ref{cor:gensmooth} 
(\ref{item:empty}) this will be the case if, and only if,
$f$ is smooth.
\end{defn}

The case, when $X$ is smooth along with $S$ is classical and at the 
basis of our results. Indeed, the study of versal maps $f:X\to S$
between smooth spaces with $T^{1}_{X/S}$ coherent over $S$ is nothing 
but the study of versal deformations of isolated
complete intersection singularities, a situation that is well 
understood. We recall the pertinent facts, with references to
Looijenga's monograph \cite{Lo4}.

\begin{theorem}
\label{thm:classical} Let $f:X\to S$ be a versal morphism between 
smooth spaces with $S$ connected. If the  critical locus
$C(f)$ is not empty, that is, $f$ itself is not smooth, then one has 
the following properties:
\begin{enumerate}
\item
\label{classical:finite}  {\sc (\cite[Thm.2.8]{Lo4})} The dimension 
of $C(f)$ equals $\dim S -1$ everywhere and $f$
restricted to $C(f)$ is birational, that is, generically one-to-one on 
each component, and finite.
\item  {\sc (\cite[Thm.4.7]{Lo4})}
$C(f)$ is locally irreducible, reduced, and determinantal, defined by 
the vanishing of the maximal minors of the Jacobi matrix
of a local embedding of $X$ into a smooth space over $S$.  In 
particular, it is locally Cohen-Macaulay. The $\calox$--module
$T^{1}_{X/S}$ is naturally a maximal Cohen-Macaulay module on $C(f)$.
\item
\label{classical:normal} {\sc (\cite[Proof of Thm.4.7]{Lo4})} The 
singular locus of $C(f)$ is either empty or equidimensional
of codimension
\begin{align*}
\codim_{C(f)}\Sing C(f) = \dim X-\dim S +3 \ge 3\,.
\end{align*} In particular, $C(f)$ is normal and normalizes the 
(reduced) discriminant $f(C(f))$.
\item  {\sc (\cite[4.C]{Lo4})} The critical locus admits a 
``Springer-like'' desingularization, given by the projection onto
the first factor of
\begin{align*} {\widetilde {C(f)}} =\{(x,H)\in C(f)\times 
\bbbp(T_{f(x)}S)\mid {\rm Im}(\jac_{x}(f))\subseteq H\subset T_{f(x)}S\}\,,
\end{align*} where we identify points $H\in  \bbbp(T_{f(x)}S)$ with 
hyperplanes in the tangent space $T_{f(x)}S$ of $S$ at
$f(x)$. The same desingularization can be obtained as the Nash 
blow-up  (or ``development'') of the discriminant, the closure
of the image of the Gau\ss\ map of the hypersurface $f(C(f))\subset S$.
\item
\label{classical:free}  {\sc (\cite[6.D]{Lo4})} The equivalent 
conditions of \ref{lem:depthT1} are satisfied and so
\ref{cor:gensmooth} applies. The $\calos$--module
$f_{*}T^{1}_{X/S}$ is canonically isomorphic to the Jacobi ideal of 
the homological discriminant that endows $f(C(f))$ with its reduced structure.
The discriminant is in particular a free divisor.
\end{enumerate}

\end{theorem}

Now we are able to formulate our main results.

\begin{theorem}
\label{thm:main1}  Let $f:X\to S$ be a versal morphism of analytic 
germs with $S$ smooth.  If $\codim_{S}f(\Sing X)\ge 2$ and
if the submodule $\call\subseteq T^{0}_{S}$  of liftable vector 
fields is a free $\calos$--module, then the discriminant of
$f$  is a {\em free divisor\/}, the liftable vector fields coincide
with the logarithmic vector fields along the discriminant, and,
as an $\calos$--module, $f_{*}T^{1}_{X/S}$ is isomorphic to the Jacobi ideal 
of the discriminant.

If even $\codim_{S}f(\Sing X)\ge 3$, then the algebra 
$\End_{\calos}(f_{*}T^{1}_{X/S})$ is the {\em normalization\/} both of the
critical locus and the discriminant of $f$, unless $f$ is smooth.
\end{theorem}

\begin{proof} The assumptions on $f(\Sing X)$ and $\call$ allow to apply \ref{cor:gensmooth}. If $f$ itself is smooth, then the discriminant is empty 
and the assertion is vacuously true.

Else, $f_{*}T^{1}_{X/S}$ is supported on a proper hypersurface in $S$ by
\ref{cor:gensmooth} and so $f':=f|_{X\setminus f^{-1}f(\Sing X)}\to 
S\setminus f(\Sing X)$ is a flat morphism between smooth
spaces that  is itself not smooth.  It follows then from 
\ref{thm:classical}(\ref{classical:free}) that $f_{*}T^{1}_{X/S}$  is
generically free of rank one on each component of its support. This 
implies that the defining equation $\Delta$ of the
discriminant of $f$  is reduced and the result follows from Saito's 
criterion \cite{Sai}, identifying at the same time $\call$ with
$\cDer_{S}(-\log\Delta)$ and $f_{*}T^{1}_{X/S}$ with the Jacobi ideal
of $\Delta$, in view of \ref{sit:Aleks}.

For the final assertion, let $ E:=\End_{\calos}(f_{*}T^{1}_{X/S})$. As $f_{*}T^{1}_{X/S}$ is coherent, isomorphic to the Jacobi ideal of $\Delta$, the $\calos$-algebra $E$ is again analytic and the structure morphism $\calo_\Delta \to E$ is finite and generically an isomorphism. Moreover, this morphism of analytic algebras factors through $\calo_{C(f)}$. The assumption on the codimension of $f(\Sing X)$ guarantees that $\calo_{C(f)} \to E$ is an isomorphism in codimension 2. There, however, $\calo_{C(f)}$ is normal by \ref{thm:classical}(\ref{classical:normal}), and so $\calo_{C(f)} \cong E$. In particular, $E$ satisfies Serre's condition $R_1$. On the other hand, $f_{*}T^{1}_{X/S}$, being maximal Cohen-Macaulay on the hypersurface $\Delta$, is a reflexive $\calos$-module and so $E$ satisfies Serre's condition $S_2$. Taken together, $E$ is normal, thus the normalization of both $\calo_\Delta$ and $\calo_{C(f)}$.
\end{proof}

\begin{remark}
Let $X$ be any singularity whose local ring satisfies Serre's condition $S_{2}$,
for example, $X$ a hypersurface, or, more generally, Cohen-Macaulay.
With $J =\jac_{X}$ the Jacobi ideal of $X$ and $J^{-1}$ its $\calox$--dual,
the endomorphism ring $E:=\End_{X}(J^{-1})$ is again an analytic algebra that
sits between $\calox$ and its normalization. 

A fundamental fact, established
by Vasconcelos \cite{Vas1} in the affine case, and, in slightly different form, earlier by
Grauert-Remmert \cite{GR1, GR2}; see also \cite{deJ}; in the analytic case, says that $E$ 
coincides with $\calox$ if, and only if, $X$ is already normal.
As $E$ inherits Serre's property $S_{2}$, this yields evidently
an algorithm for normalization. See also \cite{Vas3} for a detailed
discussion.

If the support $\Delta$ of $T^{1}_{X/S}$ is a free divisor, as in the situation of 
Theorem  \ref{thm:main1}, then it is isomorphic to the Jacobian ideal
$J$ of $\Delta$, and its $\calo_{\Delta}$--dual $J^{-1}$ is again a maximal 
Cohen-Macaulay module on $\Delta$. In particular, transposition defines an
isomorphism of analytic algebras, $\End_{\Delta}(J)\xto{\cong} \End_{\Delta}(J^{-1})$.

The second half of \ref{thm:main1} thus says that, for 
$\codim_{S}f(\Sing X)\ge 3$, the discriminant of $f$ is ``almost normal'', 
in the sense that a single step in the algorithm suffices to normalize it.
\end{remark}

To apply now the main result, we need thus criteria that guarantee freeness 
of the module of liftable vector fields and allow to bound
the singular locus of $X$. The case of a finite flat map is easy to 
analyse as we now show.

\section{Discriminants of Finite Flat Maps}
\begin{theorem} Let $f:X\to S$ be a finite, flat map of complex 
spaces with $X$ {\em normal\/} and $S$ {\em smooth}. The
discriminant of $f$ is a (non empty) free divisor at each point $s\in 
S$,  at which $f_{*}\calt^{0}_{X}$ is {\em locally
free\/}  and $f$ is {\em versal} (and ramified).
\end{theorem}

\begin{proof} To begin with, we repeat some of the pertinent 
arguments from above in the  context of sheaves of
$\calox$--modules. Applying
$\cHom_{X}(?,\calox)$ to the Zariski-Jacobi sequence for $f$  yields 
the exact sequence of $\calox$-modules
$$ 0\to \calt^{0}_{X/S}\to \calt^{0}_{X}\to
\calt^{0}_{S}(\calox)\to \calt^{1}_{X/S}\to \calt^{1}_{X}\to 0\,,
$$ where the zero at the right end is due to the smoothness of $S$. 
Moreover, $f$ being flat onto the smooth space
$S$, the space $X$ is locally Cohen-Macaulay. As $X$ is normal, thus 
reduced, $f$ is necessarily generically
\'etale by generic smoothness, and that forces
$\calt^{0}_{X/S} = 0$.

At each point $s\in S$, at which $f$ is versal, already the 
Kodaira-Spencer map $\calt^{0}_{S}\to \calt^{1}_{X/S}$ is
surjective, thus so is, as in the proof of \ref{cor:rigid},  a fortiori the map
$\calt^{0}_{S}(\calox)\to \calt^{1}_{X/S}$. In all, the (direct image 
of the) above exact sequence reduces at each versal
point to
$$ 0\to f_{*}\calt^{0}_{X}\to f_{*}\calt^{0}_{S}(\calox)\to 
f_{*}\calt^{1}_{X/S}\to 0\,.
$$ Now $f_{*}\calt^{0}_{S}(\calox) = \calt^{0}_{S}(f_{*}\calox)$ is 
locally free as
$\calos$-module, because $S$ is smooth and $X$ is locally Cohen-Macaulay. 
Moreover, $f_{*}\calt^{1}_{X/S}$ is of codimension
at least one as $f$ is generically
\'etale, and it is of projective $\calos$--dimension at most one at
$s$ if and only if $f_{*}\calt^{0}_{X} = f_{*}\Theta_{X}$ is free at 
$s$.  If that condition is satisfied, then the sheaf of
liftable vector fields
$\call\subseteq \calt^{0}_{S}$ is locally free at $s$ by 
\ref{lem:depthT1}. As $X$ is normal, $\Sing X$ is of codimension at
least two in $X$ and so is $f(\Sing X)$ in $S$ as $f$ is finite and 
flat. In all, Theorem \ref{thm:main1} applies.
\end{proof}

\begin{remark} The proof shows that conversely, if $f$ is versal at 
$s\in S$, then the  homological discriminant is a free
divisor at $s$ if, and only if, $\calt^{0}_{X}$ is a  maximal 
Cohen-Macaulay $\calox$--module at each point $x\in X$ over $s$,
equivalently, if $f_{*}\calt^{0}_{X}$ is locally free at $s$.
\end{remark} To produce free divisors from finite flat maps, one may 
thus start with some unobstructed Artinian scheme (or
space)
$X_{0}$ and consider a versal deformation $X\to S$. The only 
conditions left to be satisfied are then that $X$ is normal with
$f_{*}\calt^{0}_{X}$ locally free, equivalently,
$\calt^{0}_{X}$ a maximal Cohen-Macaulay $\calox$-module. Although it is 
rare that the module of vector fields is maximal
Cohen-Macaulay for a normal Cohen-Macaulay singularity, it happens in 
the following case.

\begin{prop} Let $X_{0}$ be an artinian space that is (algebraically) 
linked to a complete intersection. One has then the
following facts:
\begin{enumerate}
\item The deformation theory of $X_{0}$ is unobstructed, thus, any 
versal deformation $f:X\to S$ of $X_{0}$ is a finite flat
map onto a smooth base $S$.

\item The total space $X$ of any versal deformation of
$X_{0}$ is rigid, Cohen-Macaulay and nonsingular in codimension three.

\item The $\calox$-module $\calt^{0}_{X}$ is maximal Cohen-Macaulay.
\end{enumerate}
\end{prop}

\begin{proof} This result, except for the statement on the singular 
locus of $X$, was first established in \cite{Buc}.
Alternatively, the claims are easily obtained from the work of 
Huneke-Ulrich in \cite{HU1}.  Just using \cite[Thm 4.2]{HU1},
which gives inter alia the assertion on the singular locus of $X$, 
the remainder of the statements  follows as well from the
work of J.~Herzog in \cite{Her}.
\end{proof}

As an immediate consequence we have the following result.

\begin{theorem} Let $X_{0}$ be an artinian space that is 
(algebraically) linked to a complete intersection. For every versal
deformation $f:X\to S$ of $X_{0}$, the cohomological discriminant is 
a free divisor and $E
=\End_{S}(f_{*}\calt^{1}_{X/S})$ normalizes both the discriminant 
and the critical locus of $f$.\qed
\end{theorem}

An interesting case to which this result applies is the Hilbert 
scheme of points on a smooth complex surface. Indeed, if $Z$
is a smooth surface, any artinian subscheme of it is Cohen-Macaulay 
of codimension 2, thus linked to a complete intersection.
Without using linkage,  Fogarty had already shown in \cite{Fog} that 
the Hilbert scheme is smooth in this case, and versality
of the Hilbert scheme at each point as a deformation of the 
associated subscheme follows from its representability, the
classical result  established by Grothendieck. In summary, we have 
hence the following application.

\begin{cor}
\label{cor:Hilbsurfaces} 
The Hilbert scheme $\calh$ of artinian subschemes of a 
smooth surface $Z$ is smooth, the universal family $f:X\to
\calh$ is versal everywhere, and $f_{*}\calt^{0}_{X}$ is locally free on
$S$. The discriminant of $f$ is a free divisor at each 
point and the endomorphism ring of its Jacobi ideal
normalizes both discriminant and critical locus.
\end{cor}

\begin{remark}
\label{rem:Hilbcurves}
The preceding corollary holds as well, but is simpler, for the Hilbert scheme of 
artinian subschemes of a smooth {\em curve\/}. Indeed, an artinian 
subscheme of a smooth curve is a disjoint union of irreducible
zero-dimensional schemes that are then necessarily simple of type $A_{n}$,
for appropriate $n$, as the local analytic rings of the curve are just
power series rings in one variable. The Hilbert scheme at the corresponding point
is the product of versal deformations of the individual singularities. In particular, the
occurring singularities are hypersurface singularities to which \ref{thm:classical} applies.
Arnol'd, in \cite{Ar1, Ar2}, investigated the associated discriminants of those singularities, establishing that they form free divisors and giving explicit closed (local) equations.
We review and extend some of that work below.
\end{remark}

\section{The Case of Curves and Surfaces}
{\em In this section, $f:X\to S$ will always denote a {\em flat\/} morphism  of
analytic germs with $S$ {\em smooth}.}

To apply our main result \ref{thm:main1} to versal deformations of 
curve or surface singularities, we first investigate
alternative characterizations of the freeness of the module of 
liftable vector fields. The key idea, already exploited in
\cite{vSt}, is to investigate the depth of  $T^{0}_{X/S}$, the module 
of ``vertical'' vector fields with respect to $f$. To
this end, one uses the following two results, the first of which 
simply recalls the behaviour of depth in short exact
sequences.

\begin{lemma}
\label{lemma:depth} If $X$ is rigid, Cohen-Macaulay, and if the 
module of vector fields on $X$ satisfies
$$\depth_{X}T^{0}_{X}\ge\min\{\dim X, \depth_{X} T^{1}_{X/S}+2\}\,,
$$
then the depth of the module of vertical vector fields is given by
\begin{align*}
\depth_{X} T^{0}_{X/S} = \min\{\dim X,\depth_{X} T^{1}_{X/S}+2\}\,.
\end{align*}
\end{lemma}

\begin{proof} As $X$ is rigid, the initial segment of the 
Zariski-Jacobi sequence for the tangent cohomology of $f$ takes the
following form; see \ref{sit:basic}(\ref{diag:ZJKSseq}) and \ref{cor:rigid};
\begin{align*} 0\to T^{0}_{X/S}\to T^{0}_{X}\to T^{0}_{S}(\calox)\to 
T^{1}_{X/S}\to 0\,.
\end{align*}
Now use the fact that
$T^{0}_{S}(\calox)\cong T^{0}_{S}\otimes_{\calos}\calox$ is a free
$\calox$--module, thus of maximal depth as $\calox$--module, and the 
mentioned property of depth in short exact sequences.
\end{proof}

\begin{lemma}
\label{lemma:isos} Assume $X$ is Cohen-Macaulay with dualizing module 
$\omega_{X}$. If the restriction of $f$ to its critical
locus is finite but not dominant, then one has, with $d=\dim X - \dim 
S$, natural isomorphisms of $\calox$--modules
\begin{align*} T^{0}_{X/S} \cong \Hom_{X}(\Omega^{1}_{X/S},\calox)\cong
\Hom_{X}(\omega_{X},(\Omega^{d-1}_{X/S})^{**})\,,
\end{align*} where $(\ )^{*}$ denotes the $\calox$--dual.
\end{lemma}

\begin{proof} The first isomorphism can be taken as the definition of 
$T^{0}_{X/S}$. The $\calox$--module
$\Omega^{i}_{X/S}$ is locally free outside of $C(f)$, for any $i$, 
and $\Omega^{d}_{X/S}$ coincides there with
$\omega_{X/S}\cong\omega_{X}$, the last isomorphism due to smoothness 
of the germ $S$.

If the relative dimension $d:=\dim X-\dim S$ of $f$ equals zero, that 
is, $f$ finite, then $\Omega^{1}_{X/S}$ is a torsion
$\calox$--module because $C(f)$ is a proper closed subspace of $X$ by 
assumption. In this case, the claimed isomorphism is one
of zero modules.

If the relative dimension $d$ is at least $1$, the assumptions imply 
that the codimension of the critical locus of $f$ in $X$
is at least $2$; see, for example, \cite[Thm.2.5]{Lo4}. With 
$j:X\setminus C(f)\into X$ the open embedding of the complement
of the critical locus into $X$, one has then natural isomorphisms
\begin{align*} (\Omega^{d}_{X/S})^{**}\cong 
j_{*}j^{*}\Omega^{d}_{X/S}\cong \omega_{X/S}\cong \omega_{X}
\end{align*} and the $\calox$--module homomorphisms
\begin{align*}
  \Hom_{X}(\Omega^{1}_{X/S},\calox)\xto{\alpha}
\Hom_{X}(\Omega^{d}_{X/S},\Omega^{d-1}_{X/S})\xto{\beta}\Hom_{X}((\Omega^{d}_{X/S})^{**},(\Omega^{d-1}_{X/S})^{**})\,,
\end{align*} with
\begin{align*}
\alpha(\vp)(x_{0}dx_{1}\wedge\cdots\wedge dx_{d})&:= 
\sum_{i=1}^{d}(-1)^{i-1}x_{0}\vp(dx_{i}) 
dx_{1}\wedge\cdots\wedge{\widehat
{dx_{i}}}\wedge\cdots\wedge dx_{d}
\end{align*} and $\beta(\psi) := \psi^{**}$, are isomorphisms outside 
of $C(f)$, thus outside a closed subset of codimension
at least $2$. As the first and last term are reflexive 
$\calox$--modules, the composition $\beta\alpha$ is necessarily an
isomorphism of $\calox$--modules.
\end{proof}

After these preliminaries, we can now formulate the following result 
for versal families of  {\em smoothable, reduced curve
singularities\/}, extending the material in \cite{vSt} that motivates 
it. The assumptions on the fibres mean precisely that
the critical locus is finite but not dominant over the base.

\begin{theorem}
\label{thm:curves} Assume $f:X\to S$ is versal with $S$ smooth, 
$X$ Cohen-Macaulay, $\dim X-\dim S=1$, and the restriction of $f$ to its critical
locus finite but not dominant. 

If $T^{0}_{X}$ and $\omega^{*}_{X}$ 
are maximal Cohen-Macaulay $\calox$--modules, then the
discriminant of $f$ is a free divisor.
\end{theorem}
\begin{proof} The assumptions first ensure that \ref{lemma:depth} 
applies, as here $\depth T^{0}_{X}=\dim X = \dim S+1$, and
$X$ is necessarily rigid by \ref{cor:rigid}. In view of 
\ref{lemma:isos}, one has
$T^{0}_{X/S}\cong \omega^{*}_{X}$, and so by \ref{lemma:depth} the 
depth of $T^{1}_{X/S}$ is at least $\dim X-2=\dim S -1$.
As the support of $T^{1}_{X/S}$ is contained in the critical locus, 
the  assumption on that locus yields $\depth T^{1}_{X/S} =
\dim S -1$.  The claim follows from \ref{thm:main1} in view of 
\ref{lem:depthT1}.
\end{proof}

As a first application we regain the main result from \cite{vSt}:

\begin{example}
\label{ex:codim2} Let $C\subset \bbbc^{3}$ be a reduced space curve 
singularity. As first shown by Schaps, \cite{Schaps}, the
deformation theory of such singularities is unobstructed, whence any 
versal deformation of it is of the form $f:X\to S$ with
smooth base $S$ and $\dim X-\dim S = 1$.  Furthermore, the total 
space $X$ is determinantal, given by the vanishing of the
maximal minors of a generic $n\times(n+1)$ matrix for a suitable $n$. 
In particular, the singular locus of $X$ is of
codimension $4$ in $X$, if not empty.  As further the singularity of 
$C$ is isolated, $f$ is finite but not dominant when
restricted to its critical  locus. That $T^{0}_{X}$ is then maximal 
Cohen-Macaulay was already established  in \cite{Her}, and
that $\omega^{*}_{X}$ is maximal Cohen-Macaulay as well is  contained 
in the classification of maximal Cohen-Macaulay modules
of rank one  on determinantal varieties, as one finds it in 
\cite{BVe}.  That the discriminant in this case is a free divisor
is the main result in \cite{vSt},  and \cite{MvSt} yields interesting 
explicit examples of free divisors so obtained.  The
additional information we obtain here through \ref{thm:main1} is that 
$\End_{\calos}(T^{1}_{X/S})$ already normalizes critical
locus and discriminant,  as $\codim_{S}f(\Sing X)\ge \codim_{X}\Sing 
X -(\dim X-\dim S)\ge 4 -1=3$.
\end{example}

A second large class of curve singularities that give rise to free 
divisors as discriminants in their versal deformations is
provided by the following application. The particular case of 
Gorenstein curve singularities in $\bbbc^{4}$ was already
mentioned in \cite{vSt}.

\begin{prop}
\label{prop:Gcurves} If $C$ is a Gorenstein curve singularity that is 
smoothable and satisfies $T^{2}_{C}=0$, then the
discriminant in the base space of any versal deformation $f:X\to S$ 
is a free divisor. If $C$ can further be deformed into a
non-smooth isolated complete intersection singularity, then the 
endomorphism ring of the $\calos$--module $T^{1}_{X/S}$
normalizes again discriminant and critical locus.

All the foregoing assumptions are satisfied for reduced Gorenstein 
curve singularities that are algebraically linked to a
complete intersection singularity, in particular, for reduced 
Gorenstein curve singularities in $\bbbc^{4}$.
\end{prop}

\begin{proof} If a reduced Gorenstein curve singularity $C$ is 
special fibre of a flat morphism $f:X\to S$ between complex
germs with $S$ smooth, then $X$ is Gorenstein as well and so 
$\omega_{X}\cong \calox\cong \omega^{*}_{X}$ is automatically a
maximal Cohen-Macaulay module. Moreover, if $T^{2}_{C}=0$, then any 
versal deformation of $C$ is of the form $f:X\to S$ with
$S$ smooth as the deformation theory of $C$ is unobstructed. In 
addition, however, $T^{0}_{X}$ is then maximal Cohen-Macaulay
over $\calox$ as follows from \cite[Satz 2.3 and remark before 
Satz 1.6]{Her}. Thus, \ref{thm:curves} applies to yield that
any unobstructed and reduced Gorenstein curve singularity exhibits a 
free divisor as discriminant in any versal deformation.
If $C$ deforms into a non-smooth isolated complete intersection 
singularity, then
$X$ is generically smooth along the critical locus of $f$ by 
\ref{thm:classical}, whence the result on normalization of
discriminant and critical locus.

Reduced Gorenstein curve singularities that are algebraically linked 
to a complete intersection necessarily satisfy
$T^{2}_{C}=0$ , as can be deduced from
\cite[Satz 1.4]{Her} in conjunction with \cite[Thm 4.2.12]{Vas2}, a 
published account of a result in \cite{Buc}. Any versal
deformation of such a singularity has total space that is nonsingular 
in codimension $6$ by \cite{HU1}. In particular, such
singularities are smoothable and the total space is generically 
nonsingular along the critical locus, whence all the
assumptions of the first part are satisfied.
\end{proof}

We now turn to deformations of surface singularities. It is tempting 
to hope  that in analogy to \ref{ex:codim2} any isolated
Cohen-Macaulay singularity of codimension
$2$ will produce a free divisor as discriminant in its versal 
deformation. However, that is not the case, as the following
example shows.

\begin{example} The cone over the rational normal curve in 
$\bbbp^{3}$ provides at its vertex an isolated  two-dimensional
Cohen-Macaulay singularity in codimension two, to which the results 
of \cite{Schaps} mentioned in \ref{ex:codim2} apply
mutatis mutandis. In particular, the module of vector fields on the total space of
a versal deformation is Cohen-Macaulay. 
The structure of the semi-universal deformation of 
this singularity was determined by Pinkham in
\cite{Pi1}. He showed that the dimension of its smooth base 
equals $2$, but that the original singularity constitutes the
only singular fibre. Accordingly, the discriminant consists of a 
single point and is not even a divisor. The culprit is clearly $T^{0}_{X/S}$;
it has only depth $2$.
\end{example}

In light of this example, it is perhaps somewhat surprising that the results of
\ref{prop:Gcurves} for Gorenstein curve singularities indeed extend 
to a significant class of Gorenstein surface singularities. 
To formulate it, recall that any germ $X_{0}$ of a 
complex singularity can be embedded into a smooth germ, say,
$X_{0}\subseteq Z$, and that the  corresponding {\em conormal 
module\/} is the $\calo_{X_{0}}$--module $\cali/\cali^{2}$, with
$\cali\subseteq \caloz$ the ideal defining the embedding. Conormal 
modules pertaining to different embeddings into smooth
germs are stably isomorphic as $\calox$--modules, whence a property 
such as reflexivity is shared by all such conormal modules.

\begin{theorem} \label{thm:Gorenstein3}
Let $X_{0}$ be the germ of an isolated Gorenstein 
surface singularity whose conormal module is {\em
reflexive}.  If $X_{0}$ is smoothable, then each versal deformation 
$f:X\to S$ satisfies the assumptions of \ref{thm:main1},
thus, the discriminant of $f$ is a free divisor. If, moreover, $X$ is 
generically smooth along the critical locus, then
$\End_{\calos}(T^{1}_{X/S})$ normalizes critical locus and discriminant of $f$

All the foregoing assumptions are satisfied for isolated Gorenstein 
surface singularities that are linked to complete
intersections, in particular for isolated Gorenstein surface 
singularities in $\bbbc^{5}$.
\end{theorem}

\begin{proof} As $X_{0}$ is reduced Cohen-Macaulay of dimension two, 
reflexivity is the same as maximal depth for a coherent
module. In other words, the conormal module is maximal Cohen-Macaulay 
as $\calo_{X_{0}}$--module by assumption. It follows
then from \cite{Her,Wa} that the reduced Gorenstein singularity $X_{0}$ 
satisfies $T^{2}_{X_{0}}= 0$, whence the deformation
theory of such a surface singularity is unobstructed and so any 
versal deformation $f:X\to S$ has a smooth base $S$, and that
furthermore $T^{0}_{X}$ is maximal Cohen-Macaulay as 
$\calox$--module, $X$ being rigid and again Gorenstein. To apply
\ref{thm:main1},  it thus remains to verify that $T^{0}_{X/S}$ is of 
depth at least $\dim S + 1 = \dim X -1$.  To this end, we
use the isomorphism in
\ref{lemma:isos} that reduces here to $T^{0}_{X/S}\cong 
(\Omega^{1}_{X/S})^{**}$, the reflexive hull of the module of relative
differential forms, as $d=\dim X-\dim S =2$ and  as $\omega_{X}\cong 
\calox$ because $X$ is Gorenstein.

We establish that $(\Omega^{1}_{X/S})^{**}$ is of depth at least 
$\dim S + 1 = \dim X -1$  in two steps, showing first that
$\Omega^{1}_{X/S}$ itself has the desired depth and then that this 
module is reflexive.

As concerns the depth of $\Omega^{1}_{X/S}$, first note that we may 
lift any embedding $X_{0}\subseteq Z$ of the original germ
into a smooth one to an embedding $X\subseteq Z\times S$ such that 
$f$ factors into this embedding followed by the projection
onto $S$ in the second factor. Using yet again the results from 
\cite{Her}, the conormal
$\calox$--module $\calj/\calj^{2}$ with respect to the ideal 
$\calj\subseteq \calo_{Z\times S}$ defining the embedding of $X$
is  a maximal Cohen-Macaulay $\calox$--module. As $X$ is reduced 
along with $X_{0}$, this implies that the Zariski-Jacobi
sequence associated to the embedding $X\subseteq Z\times S$ over $S$ 
is {\em exact at the left\/}, that is,
\begin{align}
\label{eq:ZJexact} 0\to \calj/\calj^{2}\xto{j} \Omega^{1}_{Z\times 
S/S}\otimes_{\caloz}\calox\to \Omega^{1}_{X/S}\to 0
\end{align} is a short exact sequence. Indeed, outside the critical 
locus of $f$  this sequence is even split exact, whence
the kernel of $j$ is a torsion $\calox$--submodule of the maximal 
Cohen-Macaulay module $\calj/\calj^{2}$, thus is zero. Now
$\Omega^{1}_{Z\times S/S}\otimes_{\caloz}\calox$ is a free 
$\calox$--module, whence this short exact sequence presents
$\Omega^{1}_{X/S}$ as the cokernel of a monomorphism between maximal 
Cohen-Macaulay modules. Its depth is thus at most one
less than the dimension of $X$.

It remains to establish reflexivity of $\Omega^{1}_{X/S}$. To this 
end we use a general criterion due to M.Auslander: Let
\begin{align*}
\xymatrix{ F_{2}\ar[r]^{\beta}&F_{1}\ar[r]^{\alpha}&F_{0}\ar[r]& 
\Omega^{1}_{X/S}\ar[r]&0 }
\end{align*} be the beginning of a free $\calox$--resolution of $ 
\Omega^{1}_{X/S}$ and set
\begin{align*} Tr:= Tr( \Omega^{1}_{X/S}):= \coker \alpha^{*}\,,
\end{align*} the {\em Auslander transpose\/} of  $\Omega^{1}_{X/S}$. 
One has then an exact sequence
\begin{align*} 0\to \Ext^{1}_{\calox}(Tr,\calox)\to 
\Omega^{1}_{X/S}\to (\Omega^{1}_{X/S})^{**}
\to \Ext^{2}_{\calox}(Tr,\calox)\to 0\,,
\end{align*} with the morphism in the middle the canonical map into 
the reflexive hull. Thus, $ \Omega^{1}_{X/S}$ is reflexive
if, and only if,
$\Ext^{i}_{\calox}(Tr,\calox)=0$ for $i=1,2$. As $\Omega^{1}_{X/S}$ 
is cokernel of a monomorphism of reflexive modules in view
of the short exact sequence  (\ref{eq:ZJexact}) above, it follows 
immediately, say,  from the snake lemma, that
$\Omega^{1}_{X/S}$ embeds into its reflexive hull,  thus, 
$\Ext^{1}_{\calox}(Tr,\calox)=0$.  To establish vanishing of
$\Ext^{2}_{\calox}(Tr,\calox)$, we may choose 
$F_{0}=\Omega^{1}_{Z\times S/S}\otimes_{\caloz}\calox$  and identify
$\coker\beta$ with $\calj/\calj^{2}$, thus
$\ker\beta^{*}\cong \caln_{X/Z\times S}:= 
\Hom_{\calox}(\calj/\calj^{2},\calox)$,  the {\em normal module\/} of 
the embedding
of $X$ into $Z\times S$.  Dualizing the displayed initial segment of 
the free resolution into $\calox$  yields then short
exact sequences
  \begin{align*} 0\to \Image\alpha^{*}\to &\ker\beta^{*}\cong 
\caln_{X/Z\times S}
\to \Ext^{1}_{X}(\Omega^{1}_{X/S},\calox)\to 0\,,\\ 0\to 
\Image\alpha^{*}\to &F^{*}_{1}\to Tr\to 0
\end{align*} that imply
\begin{align*}
\Ext^{2}_{\calox}(Tr,\calox)&\cong \Ext^{1}_{\calox}(\Image\alpha^{*},\calox)
\intertext{as $F_{1}^{*}$ is free, and}
\Ext^{1}_{\calox}(\Image\alpha^{*},\calox)&\cong 
\Ext^{1}_{\calox}(\caln_{X/Z\times S},\calox)
\end{align*} as long as
\begin{align}
\label{eq:vanishing}
\tag{$*$}
\Ext^{i}_{\calox}(\Ext^{1}_{X}(\Omega^{1}_{X/S},\calox),\calox) =0
\end{align} for $i=1,2$. Now $X_{0}$ is an isolated smoothable 
surface singularity, whence the critical locus of $f$ is of
codimension at least $3$ in $X$. As 
$\Ext^{1}_{X}(\Omega^{1}_{X/S},\calox)$ is concentrated on $C(f)$, 
the required vanishing
in (\ref{eq:vanishing}) follows. To conclude the argument, we finally 
use that with $\calj/\calj^{2}$ also its
$(\omega_{X}=\calox)$--dual module $\caln_{X/Z\times S}$ is maximal 
Cohen-Macaulay, thus satisfies
$\Ext^{i}_{\calox}(\caln_{X/Z\times S},\calox) = 0$ for $i\neq 0$.

The final assertion on linkage follows as for curves, it simply 
exploits that we know additionally that $X$ is smooth in
codimension $6$.
\end{proof}


\section{The Classical Discriminant of a Polynomial}

\newcommand{\pdual}{\check{p}}
\newcommand{\suchthat}{\,|\,}
\newcommand{\problem}[1]{{\bf (#1)}}

We first note a classical relation between discriminants and dual varieties. For this
consider the incidence variety
\[
	I = \{ (p,H) \suchthat p \in H \} \subset \PP^n \times \PPdual^n
\]
together with the two natural projections $p \colon I \to \PP^n$ and $\pdual \colon I \to \PPdual^n$.

\begin{proposition} \label{prop:DualDiscriminant}
Let $C \subset \PP^n$ be a projective curve and 
\[
	I_C := p^{-1}(C) = \{ (p,H) \suchthat p \in C \cap H \} 
	\subset C \times \PPdual^n \subset \PP^n \times \PPdual^n
\]
the corresponding incidence variety. Then the dual variety $D \subset \PPdual^n$ of $C$ is
the discriminant of the morphism
$
	\pdual|_{I_C} \colon I_C \to \PPdual^n.
$
\end{proposition}

\begin{proof}
The fiber of $\pdual|_{I_C}$ over a point $H \in \PPdual^n$ is the intersection $C \cap H \subset \PP^n$.
It contains a point of multiplicity at least two if and only if $H$ is tangent to $C$.
\end{proof}

This gives the following well known description of the classical discriminant of a polynomial:

\begin{cor}
The discriminant of the universal polynomial
\[
	F(u,v):=s_0 u^n + \dots + s_n v^n 
\]
is isomorphic to the dual variety of the rational normal curve $\PP^1 \hookrightarrow \PP^n$
of degree $n$.
\end{cor}

\begin{proof}
Choosing coordinates $y_i$ of $\PP^n$ and dual coordinates $s_i$ of $\PPdual^n$ the incidence
variety $I$ is described by $\sum_i s_iy_i = 0$ in $\PP^n \times \PPdual^n$. If we choose
coordinates $(u:v)$ of $\PP^1$, the $d$-uple embedding is given by $(u:v) \mapsto (u^n:\dots:v^n)$.
Therefore the equation of $I_{\PP^1} \subset \PP^1 \times \PPdual^n$ is 
\[
	s_0u^n + \dots + s_n v^n = 0,
\]
the universal polynomial. The corollary follows from \ref{prop:DualDiscriminant}
\end{proof}

\begin{remark}
Notice that $\pdual_{|_C} \colon I_C\to \PPdual^{n}$ is 
the {\em universal family\/} over the {\em Hilbert scheme\/} 
$\Hilb^{n}_{\PP^1}\cong \PPdual^{n}$ of subschemes of length $n$
on $\PP^1$. 
\end{remark}

\begin{remark}
Note that transverse to the  rational normal curve\footnote{
We assume here that the binomial coefficients 
$\binom{n}{i}$ are {\em invertible}.}
\begin{align*}
\bbbp^{1}\ni (a,b)\mapsto (a^{n}: \cdots : \binom{n}{i}a^{n-i}b^{i} : \cdots : b^{n})\in \PPdual^{n}
\end{align*}
we find a semi-universal deformation of the $A_{n-1}$--singularity $(ax-by)^{n}=0$.
\end{remark}

The following homogeneous equation for the discriminant is due to Bezout (see e.g.\ \cite{GKZ}).
\begin{proposition}
The coefficients $s_{ij}$ of the generating function
\begin{align*}
\sum_{i,j=1}^{n-1}s_{ij}x^{n-1-i}y^{n-1-j} := 
\frac{F_{v}(x,1)F_{u}(y,1)-F_{v}(y,1)F_{u}(x,1)}{y-x}
\end{align*}
are homogeneous quadratic polynomials in $\bbbz[s_{0},...,s_{n}]$ 
and the {\em Bezout determinant\/} of the symmetric 
$((n-1)\times(n-1))$--matrix $\bbbb := (s_{ij})_{i,j=1,...,n-1}$,
\begin{align*}
B(s_{0},...,s_{n}):= \det \bbbb
\end{align*}
yields an equation of the discriminant in $\PPdual^{n}$. It is 
homogeneous of degree $2n-2$,
as well as weighted homogeneous of degree $n(n-1)$ with 
respect to the weights $w(s_{i}) = i$.

The specialization $\bbbb|_{s_{0}=1, s_{1}=0}$ constitutes a 
{\em discriminant matrix\/} for the discriminant in the semi-universal 
deformation $f(x) = x^{n}+s_{2}x^{n-2}+\cdots + s_{n}$ of $x^{n}=0$ over 
$S=\bbbc^{n-1}_{s_{2},...,s_{n}}$.\qed
\end{proposition}

While the above yields a description of the (homogeneous) 
discriminant in the Hilbert scheme as the determinant of a square matrix 
of size $n-1$, one may also consider the following slightly less 
economical version that yields the 
discriminant of the versal deformation.

To this end, recall the Euler relation $nF(u,v) = uF_{u}(u,v) +vF_{v}(u,v)$, 
or, in its dehomogenized form, $F_{v}(x,1) = nF(x,1) - xF_{u}(x,1)$. 
If one replaces the partial derivative
$F_{v}$ by $F$ in the above expression for the Bezout form, 
one finds the following result.

\begin{theorem} \label{thm:Bezoutn}
The coefficients $s'_{ij}$ of the generating function
\begin{align*}
\sum_{i,j=1}^{n}s'_{ij}x^{n-i}y^{n-j} := \frac{F(x,1)F_{u}(y,1)-F(y,1)F_{u}(x,1)}{y-x}
\end{align*}
are again homogeneous quadratic polynomials in $\bbbz[s_{0},...,s_{n}]$ 
and the determinant of the  symmetric $(n\times n)$--matrix $\bbbb' := 
(s'_{ij})_{i,j=1,...,n}$ satisfies
\begin{align*}
\det \bbbb' = s_{0}^{2}B(s_{0},...,s_{n})\,.
\end{align*}
The matrix $\bbbb'|_{s_{0}=1}$ is a discriminant matrix for the 
discriminant in the versal deformation 
$f(x) = x^{n}+s_{1}x^{n-1}+\cdots +s_{n}$ of $x^{n}=0$ over 
$S=\bbbc^{n}_{s_{1},...,s_{n}}$.

Moreover, in the case $s_0=1$, if we split $F(x,1) = f(x) = (x-r_{1})\cdots(x-r_{n})$, then 
the entries $s'_{ij}$ satisfy
\begin{align}
\label{Weyl}
\tag{$*$}
s'_{ij} = \langle \grad_{\bf r}s_{i}, \grad_{\bf r}s_{j}\rangle = 
\sum_{k=1}^{n}\frac{\partial s_{i}}{\partial r_{k}}\cdot \frac{\partial s_{j}}{\partial r_{k}}\,.
\end{align}
\end{theorem}

\begin{proof} 
Formula (\ref{Weyl}) follows from the following elementary calculation (see also \cite[(2.4.5)~Lemma]{SYS}):
\begin{eqnarray*}
\frac{F(x,1)F_{u}(y,1)-F(y,1)F_{u}(x,1)}{y-x} & = & \frac{1}{y-x} \sum_{k=1}^n \left( -\frac{f(x)f(y)}{y-r_k} + \frac{f(x)f(y)}{x-r_k} \right)  \\
& = & \sum_{k=1}^n \frac{f(x)}{x-r_k} \frac{f(y)}{y-r_k} = \sum_{k=1}^n \frac{\partial f}{\partial r_k}(x) \frac{\partial f}{\partial r_k}(y) \\
& = & \sum_{i,j=1}^{n} \left( \sum_{k=1}^n \frac{\partial s_{i}}{\partial r_{k}}\cdot \frac{\partial s_{j}}{\partial r_{k}} \right) x^{n-i}y^{n-j}.
\end{eqnarray*}
Here we have used 
$$\frac{\partial f}{\partial r_i}(x) = - \prod_{j \neq i} (x- r_j).$$
Let $\bbbm$ be the $(n \times n)$-matrix 
$\bbbm = (\frac{\partial s_{i}}{\partial r_{j}})_{i,j=1,...,n}$. Then we have $\bbbb'|_{s_{0}=1} = \bbbm \bbbm^T$. 
Let 
$$\bbbv = \begin{pmatrix} r_1^{n-1} & \cdots & r_{n}^{n-1} \\
\vdots &  & \vdots \\
r_1 & \cdots & r_{n} \\
1  & \cdots & 1
\end{pmatrix}$$
be the (reversed) Vandermonde matrix. It has determinant
$$\det \bbbv = (-1)^{\frac{1}{2}n(n-1)} \prod_{n \geq i>j \geq 1} (r_i-r_j).$$
Then we have
$$\bbbv^T \bbbm = \begin{pmatrix} \frac{\partial f}{\partial r_1}(r_1) & \cdots & \frac{\partial
f}{\partial r_{n}}(r_1)\\
\vdots & \ddots & \vdots \\
\frac{\partial f}{\partial r_1}(r_{n}) & \cdots & \frac{\partial f}{\partial r_{n}}(r_{n})
\end{pmatrix},$$
whence
$$\det \bbbv^T \bbbm = (-1)^n \prod_{n \geq i>j \geq 1} (r_i-r_j)^2.$$
From this we obtain
$$\det \bbbb'|_{s_{0}=1} = \det \bbbm \bbbm^T = \prod_{n \geq i>j \geq 1} (r_i-r_j)^2=: \Delta.$$
By the chain rule
$${\rm grad}_{\bf s} ( \log \Delta) \bbbm \bbbm^T = {\rm grad}_{\bf r}( \log \Delta) \bbbm^T.$$
Since the components of the latter vector are symmetric polynomials in $r_1, \ldots, r_{n}$, it follows that
$${\rm grad}_{\bf s} ( \log \Delta) \bbbb'|_{s_{0}=1} \in \QQ[s_1, \ldots ,s_n]^n.$$
This shows that $\bbbb'|_{s_{0}=1}$ is a discriminant matrix for the 
discriminant in the versal deformation 
$f(x) = x^{n}+s_{1}x^{n-1}+\cdots +s_{n}$ of $x^{n}=0$ over 
$S=\bbbc^{n}_{s_{1},...,s_{n}}$.
\end{proof}

The description of the entries of the discriminant matrix in terms 
of the derivatives of the elementary symmetric functions $s_{i}$ with 
respect to the roots $r_{k}$ is precisely the form
of the discriminant matrix given by Arnol'd in \cite{Ar1,Ar2}, and its relation 
to the Bezout form can be found, at least implicitly, in \cite{SYS}. 
The form (\ref{Weyl}) of the entries of the discriminant matrix generalizes 
both to simple hypersurface singularities and to the simple elliptic surface 
singularities, in that it uses the action of the associated Coxeter or Weyl 
group and the fact that the discriminant is precisely the image of the 
union of the reflection hyperplanes under the orbit map.

\begin{remark}
Recently C.~D'Andrea and J.~V.~Chipalkatti \cite{An} proved that even the cone over the dual variety of the rational normal curve of degree $n \geq 3$ is a free divisor.
\end{remark}


\section{The dual variety of an elliptic normal curve} \label{section:EllipticNormal}

\newcommand{\AZ}{\mathbb{A}}

Let  $D_\tau$ be the dual variety of the elliptic normal curve $E_\tau \subset \PP^n$. 
In this section we find a new determinantal expression for $D_\tau \subset \PPdual^n$.

For this we turn the method of the previous section around and identify $D_\tau$ with 
the discriminant of the universal meromorphic function of degree $n+1$ with only one pole 
at $0$ on $E_\tau$. Recall that the space of such meromorphic functions is generated by
the Weierstrass $\wp$-function, its derivatives  and the constant function.

\begin{prop} \label{prop:Lambda}
Let $\tau$ be a point in the upper half plane, $E_\tau$ the corresponding elliptic curve and 
\[
	\begin{array}{ccccl} \Phi & : & E_\tau & \to & \PP^n \\ & & z & \mapsto 
	& (1:\wp(z):\wp'(z): \ldots :\wp^{(n-1)}(z))
	\end{array} 
\]
its $n$-th Weierstrass embedding. Then the dual variety of $E_\tau$ is isomorphic to the discriminant of
the universal meromorphic function 
\begin{equation}
\lambda(z)=\frac{(-1)^{n-1}}{n!} \y_{n+1} \wp^{(n-1)}(z)  \pm  \cdots - \frac{1}{2} \y_3 \wp'(z) + \y_2\wp(z) +\y_0
\label{eqLambda}
\end{equation}
\end{prop}

\begin{proof}
If we denote by $\frac{(-1)^i}{i!}s_i$ the dual coordinates of $\PPdual^n$, then the claim follows from
\ref{prop:DualDiscriminant}.
\end{proof}


Up to a constant
$u$, a meromorphic function $\lambda$ as above is determined  by its zeros $z_0,\dots,z_n$. The condition that the only pole should lie at $0$ implies that
$\sum_{j=0}^n z_j = 0$. Using the coordinates ${\bf v} = ( v_1, \ldots , v_n )$ with 
$$
v_k = \sum_{j=0}^{k-1} z_j = - \sum_{j=k}^n z_j, \quad k=1, \ldots , n,
$$
we obtain

\begin{prop}
The coefficient
$
\y_i = \y_i(\tau, u , {\bf v}) 
$
of the universal meromorphic function 
$$
\lambda(z)=\frac{(-1)^{n-1}}{n!} \y_{n+1} \wp^{(n-1)}(z)  \pm  \cdots - \frac{1}{2} \y_3 \wp'(z) + \y_2\wp(z) +\y_0
$$
is a
Jacobi form of weight $-i$ and index 1. If we set
$\y_{-1}(\tau, u,{\bf v}):= \tau$ then the algebra $J$ of Jacobi forms
is freely generated by $\y_{-1}, \y_0, \y_2, \ldots, \y_{n+1}$.
\end{prop}

\begin{proof} This follows from a theorem of K.~ÊWirthm\"uller \cite[(3.6) Theorem]{Wi} (see also
\cite[Section 1.4, in particular Theorem 1.4]{Be2}).
\end{proof}

\begin{defn}
Let $\Omega^1_J$ be the $J$-module of 1-forms  and 
$$I: \Omega^1_J \times \Omega^1_J \to J$$ 
be the symmetric bilinear form defined by
$$I(du,du)=I(d\tau,d\tau)=0, \quad I(du,d\tau)=I(d\tau,du)=1,$$
$$I(dv_i,dv_j)=\delta_{ij}, I(du,dv_i)=I(d\tau,dv_i)=0, \quad i,j
\in\{1, \ldots, n\}.$$
\end{defn}

\begin{prop} \label{prop:discLambda}
The matrix
$$B:=(I(d\y_i,d\y_j))^{i=-1,0,2, \ldots , n+1}_{j=-1,0,2, \ldots , n+1}$$
is a discriminant matrix for the discriminant of the
 universal meromorphic function 
\begin{equation}
\lambda(z)=\frac{(-1)^{n-1}}{n!} \y_{n+1} \wp^{(n-1)}(z)  \pm  \cdots - \frac{1}{2} \y_3 \wp'(z) + \y_2\wp(z) +\y_0
\end{equation}
\end{prop}

\begin{proof}
Since $\y_{-1},\y_0,\y_2,\dots, \y_{n+1}$ generate $J$ it follows from 
\cite[Formula~(5.2.5)]{Sa2}.
that $f = \det B$ is an equation for the discriminant.
By
\cite[(5.3) Assertion]{Sa2} we have
$$B \cdot (\partial_1, \ldots ,
\partial_{n+2})^t \in \Theta_J(-\log f)^{n+2}.$$
\end{proof}

To calculate $B$ explicitly one needs to express in
terms of the $\y_i$ their derivatives. Here we pursue a different approach. For a slight variant of $I$, Bertola has given a Bezout type formula expressing  the generating function for the entries of the matrix $B$ in terms of the $\y_i$  \cite[Theorem~1.5]{Be2}. The matrix $B$ for $n=2$ is computed in  \cite[Example 1.2]{Be2}. 
Retracing his steps with {\sc Maple} for $n=4$ and dividing the first row and the first column
of the resulting matrix by $-2 \pi \sqrt{-1}$ and the other rows and columns by $e^{-2 \pi 
\sqrt{-1} u}$ we obtain:

\begin{example} \label{exA}
Let $\tau \in \HH$ be a point in the upper half-plane and $E_\tau \subset \PP^4$ the
corresponding elliptic curve in its $4$-th Weierstrass embedding. Then the dual variety
of $E_\tau$ is defined by the determinant of the matrix
$$A:=\left( \begin{array}{cccccc} 0 & \y_0 &
\y_2 &
  \y_3 &  \y_4  &
\y_5 \\
  \y_0 & a_{00} & a_{02} & a_{03} & a_{04} & a_{05} \\
  \y_2 & a_{02} & a_{22} & a_{23} & a_{24} & a_{25} \\
  \y_3 & a_{03} & a_{23} & a_{33} & a_{34} & a_{35} \\
  \y_4 & a_{04} & a_{24} & a_{34} & a_{44} & a_{45} \\
  \y_5 & a_{05} & a_{25} & a_{35} & a_{45} & a_{55}
\end{array} \right)$$ where
\begin{eqnarray*} a_{00} & := & - \frac {1}{6}g_2\y_0\y_2 - \frac 
{1}{2}g_3\y_2^{2} - \frac {1}{18}g_2^{2}\y_2\y_4 + \frac
{1}{24}g_2^{2}\y_3^{2} + \frac {1}{8}g_2g_3\y_3\y_5\\ & & \mbox{} - 
\frac {1}{12}g_2g_3\y_4^{2} + \left( \frac
{1}{288}g_2^{3}+ \frac {3}{40}g_3^{2} \right)
\y_5^{2}
\\ a_{02} & := & - \frac {1}{3}g_2\y_2^{2} - g_3\y_2\y_4 + \frac 
{3}{4}g_3\y_3^{2} + \frac{1}{12} g_2^{2} \y_3\y_5 - \frac
{1}{18}g_2^{2}\y_4^{2} + \frac{9}{80} g_2 g_3 \y_5^{2}
\\ a_{03} & := & - \frac {5}{12}g_2\y_2\y_3 + \frac {1}{2}g_3\y_3\y_4 
- \frac {4}{5}g_3\y_2\y_5 - \frac {1}{36}g_2^{2}\y_4\y_5
\\ a_{04} & := &  - \frac {1}{2}g_2\y_2\y_4 + \frac 
{21}{20}g_3\y_3\y_5 - \frac {1}{2}g_3\y_4^{2} + \frac
{1}{24}g_2^{2}\y_5^{2}
\\ a_{05} & := & - \frac {7}{12}g_2\y_2\y_5 - \frac {3}{5}g_3\y_4\y_5
\end{eqnarray*}
\begin{eqnarray*} a_{22} & := & - 2\y_0\y_2 - \frac {2}{3}g_2\y_2\y_4 
+ \frac {1}{2}g_2\y_3^{2} + \frac {3}{2}g_3\y_3\y_5 -
g_3\y_4^{2} + \frac {3}{40}g_2^{2}\y_5^{2}
\\ a_{23} & := & - 3\y_0\y_3 - \frac {8}{15}g_2\y_2\y_5 + \frac 
{1}{3}g_2\y_3\y_4 - \frac {1}{2}g_3\y_4\y_5
\\ a_{24} & := & -4\y_0\y_4 + \frac {7}{10}g_2\y_3\y_5 - \frac 
{1}{3}g_2\y_4^{2} + \frac {3}{4}g_3\y_5^{2}
\\ a_{25} & := & - 5\y_0\y_5 -  \frac {2}{5}g_2\y_4\y_5
\end{eqnarray*}
\begin{eqnarray*} a_{33} & := &  - 4\y_0\y_4 + \frac {6}{5}\y_2^{2} - 
\frac {1}{6}g_2\y_3\y_5 + \frac {1}{3}g_2\y_4^{2} -
\frac {1}{2}g_3\y_5^{2}
\\ a_{34} & := & \frac {4}{5}\y_2\y_3 - 5\y_0\y_5 + \frac 
{1}{12}g_2\y_4\y_5 \\ a_{35} & := &  \frac {2}{5}\y_2\y_4  - \frac
{1}{3}g_2\y_5^{2}
\end{eqnarray*}
\begin{eqnarray*} a_{44} & := & - 2\y_2\y_4 + \frac {6}{5}\y_3^{2} + 
\frac {1}{2}g_2\y_5^{2} \\ a_{45} & := & -3\y_2\y_5 +
\frac {3}{5}\y_3\y_4
\end{eqnarray*}
\begin{eqnarray*} a_{55} & := & -2\y_3\y_5 + \frac{4}{5}\y_4^2
\end{eqnarray*} 
with $g_2$ and $g_3$ the well-known functions
\begin{eqnarray*} g_2(\tau) & = & 60 {\sum_{m,n}}' 
\frac{1}{(m+n\tau)^4}, \\ g_3(\tau) & = & 140 {\sum_{m,n}}'
\frac{1}{(m+n\tau)^6}.
\end{eqnarray*}
As usual the symbol $\sum'$ indicates summation 
over the nonzero elements of $\ZZ \times \ZZ$.
\end{example}


\section{The discriminant of the simple elliptic singularity
$\widetilde{A}_4$} \label{section:A4}

Consider the cone $(X,0)$ over an elliptic normal curve
$E_\tau \subset \PP^4$.   The ideal  
of $E_\tau \subset \PP^4$ is described by the $4 \times
4$-Pfaffians of a skew symmetric $5 \times 5$-matrix $M$. The germ $(X,0)$ is therefore a
Gorenstein surface singularity of codimension 3, namely the
simple elliptic singularity of type
$\widetilde{A}_4$.  By \ref{thm:Gorenstein3} the discriminant of the
semi-universal deformation of $(X,0)$ is a free divisor. 

In this section we will show that the matrix $A$ from
\ref{exA} is a discriminant matrix for this discriminant. More precisely, the discriminant is isomorphic to the affine cone over
the dual variety of the elliptic normal curve $\Jac_2 E_\tau$. For this we make extensive use of the fact that 
the pfaffian description of $E_\tau$ exhibits $E_\tau$ as a linear section of the Grassmannian variety $\GG(5,2)$, that the dual variety of 
$\GG(5,2)$ is again a Grassmannian of the same type \cite{Mu}, and that the deformations of $(X,0)$
can be obtained by perturbing the entries of $M$.

We start by recalling some facts about dual varieties and linear sections.
Let $\PP^n = \PP(A)$ be a projective space with $A = H^0(\sO(1))$ and 
$\PPdual^n := \PP(A^*)$ its dual space. For every linear
subspace $\PP^m = \PP(B) \subset \PP(A) = \PP^n$ we obtain  quotient 
$A \xrightarrow{\phi} B \to 0$ and consider
\[
	\Bperp = (\ker \phi)^*.
\] By construction this gives a quotient $A^* \to B_\perp \to 0$ and 
therefore a linear subspace $\PPperp^m := \PP(\Bperp)
\subset \PP(A^*) =: \PPdual^n$ which is called the orthogonal space 
of $\PP^m \subset \PP^n$.

\begin{lemma} \label{lDualSingular} Let $X \subset \PP^n$ be a smooth 
variety, $\Xdual \subset \PPdual^n$ the dual variety,
$\PP^m \subset \PP^n$ a linear subspace and $\PPperp^m \subset 
\PPdual^n$ the orthogonal space.  Assume that $Y = X \cap
\PP^m$ contains a point $y$ with $\codim T_{y,Y} < \codim X$. Then 
there exists a point $\yperp$ in $\Yperp = \Xdual \cap
\PPperp^m$ with $\codim T_{\yperp,\Yperp} < \codim \Xdual$.
\end{lemma}

\begin{proof} If $y \in Y$ is a point as above there exists a 
hyperplane in $\PP^n$ that contains $\PP^m$ and is tangent to $y
\in X$. This hyperplane represents a point $\yperp \in \Xdual$ 
because of the tangency condition. Also $\yperp$ is in
$\PPperp^m$ since it contains $\PP^m$. We have $\yperp \in \Yperp$. 
By the symmetry of the duality correspondence the
hyperplane in $\PPdual^n$ represented by $y$ is tangent to $\Xdual$ 
in $\yperp$. It also contains $\PPperp^m$ since it is an
element of
$\PP^m= (\PPperp^m)_\perp$. Therefore $\codim T_{\yperp,\Yperp} < 
\codim \Xdual$ as claimed.
\end{proof}

We now specialize to the case of a Grassmannian in its P\"ucker 
embedding. For this  let $V$ be a $5$-dimensional vector space
and $\PP^9 := \PP(\bigwedge^2 V)$ the projective space of $2$-forms 
with coordinates $v_{ij} = v_i \wedge v_j$. The
Grassmannian
$\GG := G(V,2) \subset \PP^9$ is defined by the $4 \times 4$ 
pfaffians of the generic  skew symmetric $5 \times 5$ matrix
\[
	\begin{pmatrix}
	   0 & v_{12} & v_{13} & v_{14} & v_{15} \\
	   -v_{12} & 0 & v_{23} & v_{24} & v_{25} \\
	   -v_{13} & -v_{23} & 0 & v_{34} & v_{35} \\
	   -v_{14} & -v_{24} & -v_{34} & 0 & v_{45} \\
	   -v_{15} & -v_{25} & -v_{35} & -v_{45} & 0 \\
	\end{pmatrix}
\] The dual variety of $\GG$ is again a Grassmannian of the same type
\[
	\GGdual \cong G(2,V) \subset \PP(\bigwedge^2 V^*) =: \PPdual^9.
\] We denote the coordinates of $\PPdual^9$ by $v^\ast_{ij} = v_i^* 
\wedge v_j^*$. The incidence variety
$H \subset \PP^9 \times \PPdual^9$ is defined by the equation
\[
	\sum_{1\le i < j \le 5} v_{ij}v^\ast_{ji} = 0.
\] Let now $W$ be a $5$-dimensional quotient space of $\bigwedge^2 V$,
$\PP^4 := \PP(W) \subset \PP^9$ its projectivization and
$\PPperp^4 := \PP(\Wperp) \subset \PPdual^9$ the corresponding 
orthogonal space.

\begin{proposition} If $E := \PP^4 \cap \GG$ is smooth  of dimension 
$1$ then $E \subset \PP^4$ is an  elliptic normal curve
and $\Eperp := \PPperp^4 \cap \GGdual$ is naturally isomorphic to
$\Jac_2 E$.
\end{proposition}

\begin{proof} Since $E$ is of expected codimension we have $\deg E = 
\deg \GG =5$. Adjunction shows that the arithmetic genus
of $E$ is $1$. By \ref{lDualSingular} above $\Eperp$ must also 
be smooth and of expected codimension, i.e. an elliptic
normal curve.

For the identification of $\Eperp$ with $\Jac_2 E$ consider the 
universal quotient bundle $\sQ$ on
$\GG$ and its restriction
$\sQ_E$ to $E$. As the intersection $\PP^4 \cap \GG$ is transversal,
a locally free resolution of $\sQ_E$ is obtained by 
tensoring the Koszul complex associated to $W$ with $\sQ$
\[
	0 \to \sQ(-5) \to 5 \sQ(-4) \to 10 \sQ(-3) \to 10 \sQ(-2) \to 
5 \sQ(-1) \to \sQ \to \sQ_E \to 0.
\] 
By the Theorem 
of Bott \cite{Bot}  the cohomology of $\sQ(-n)$ vanishes
for $1 \le n \le 5$. This shows that
\[
	H^0(\sQ_E) = H^0(\sQ) = V.
\] Note that $\GGdual$ corresponds to decomposable $2$-forms $v 
\wedge v'$  and such a point lies on
$\Eperp$ if it is in the kernel of the map
\[
    \bigwedge^2 V = \bigwedge^2 H^0(\sQ_E) \to H^0\Bigl(\bigwedge^2 
\sQ_E\Bigr) =
    H^0\bigl(\sO_E(1)\bigr) = W.
\] In particular $\sQ_E$ can not have any subbundles $\sL$ of degree 
$\deg \sL \ge 3$ since in this case
$h^0(\sL) \ge 3$ by Riemann-Roch and  $\bigwedge^2 H^0(\sL) \subset 
\bigwedge^2 V$ spans a projective space of dimension ${
\deg \sL \choose 2}-1\ge 2$ that is contained in $\Eperp$. This 
contradicts $\dim \Eperp = 1$.

 From $\deg \sQ_E = 5$ it follows that $\sQ_E$ is stable. By Atiyah's 
classification of stable vector bundles on elliptic
curves \cite{Ati}, $\sQ_E$ is the unique irreducible rank $2$ vector 
bundle with determinant $\det \sQ_E = \sO_E(H)$, where $H$
is a hyperplane section of $E \subset \PP^4$.

Let now $\sL \in \Jac_2 E$ be a line bundle of degree $2$. Then there 
exists a unique nontrivial extension
\[
	0 \to \sL \to \sF \to \sO_E(H) \otimes \sL^{-1} \to 0.
\] By Atiyah's classification we must have $\sF \cong \sQ_E$. Taking 
cohomology we obtain a two dimensional subspace $H^0(\sL)
\subset V$ and a one dimensional subspace
\[
	\bigwedge^2 H^0(\sL) \subset \bigwedge^2 V
\] that is mapped to $H^0(\bigwedge^2 \sL) = 0$ by the map 
$\bigwedge^2 V \to W$. It therefore represents a projective point
on $\Eperp$.

If on the other hand $v \wedge v'$ represents a point on $\Eperp$ 
then $v \wedge v'$ is mapped to zero in $W$.  This means
that $v$ and $v'$ are dependent on $E$ and that the image of
\[
	\sO \oplus \sO \xrightarrow{(v,v')} \sQ_E
\] is a line bundle $\sL$ on $E$ with at least $2$ sections. 
Therefore the degree of $\sL$ is at least $2$. Since a subbundle
of $\sQ_E$ has degree at most $2$ we have $\sL \in \Jac_2 E$.
\end{proof}

Let now $W \to U$ be a $4$-dimensional quotient space of $W$ and 
$\PP^3 := \PP(U) \subset \PP^4$ the corresponding
hypersurface. The associated orthogonal space
\[
	\PP^5 := \PPperp^3 := \PP(\Uperp) \subset \PPdual^9
\] contains $\PPperp^4$. We set $Z := E \cap \PP^3 = \GG \cap \PP^3$ and
$Y := \GGdual \cap \PP^5$. If $E = \GG \cap \PP^4$ is a smooth 
elliptic curve, then $Z$ is a scheme of length $5$ on $E$ and
$Y$ is a possibly singular Del Pezzo surface of degree $5$ that 
contains $\Eperp$.

\begin{proposition} \label{pBothSingular} 
In the situation just described $Z$ contains a multiple point if, and only if, $Y$ is singular.
\end{proposition}

\begin{proof} This follows from \ref{lDualSingular}.
\end{proof}

\begin{proposition} \label{pDualDiscriminant} Let $\PPdualblowup^9$ 
be the blowup of $\PPdual^9$ in $\PPperp^4$ and
$\GGdualblowup$ the strict transform of $\GGdual$. Then there exists 
a natural map
$\GGdual \to \PPdual^4:=\PP(W^*)$ of fiber dimension $2$ whose 
discriminant is equal to the  dual variety $\Edual$ of
$E\subset \PP^4$.
\end{proposition}

\begin{proof} The exact sequence
\[
	0 \to W^* \to \bigwedge^2 V^\ast \to \Wperp \to 0
\] yields under projection from $\PPperp^4 = \PP(\Wperp)$ a morphism
$$\xymatrix{
	\GGdualblowup \ar[r] \ar[dr] & \PPdualblowup^9 \ar[d]^{\pi} \\
	&\PP(W^*).
	}
$$ The fiber over a point $\PP^3 \in
\PPdual^4=\PP(W^*)$ is then $\PPperp^3 = \PP^5$ and  $\PP^5 \cap 
\GGdualblowup = \PP^5 \cap \GGdual = Y$ is a singular surface
  if and only if $Z = E \cap \PP^3$ contains a double point by 
\ref{pBothSingular}. This is the case if and only if
$\PP^3$ is tangent to $E$ so $\PP^3$ is a point of the dual variety 
$\Edual$ of $E$.
\end{proof}

If we choose a splitting $\bigwedge^2 V^\ast \cong \Wperp \oplus W^*$ 
we can write every decomposable form as
\[
	v^*_i \wedge v^*_j = w^\perp_{ij} + w^*_{ij}.
\]

\begin{cor} \label{CorHC} There exists a flat deformation of the projective closure 
of the cone over $\Eperp$
$$\xymatrix{
	X \ar[r] \ar[d]^{\pi} & \PP^5 \times W^* \ar[dl]\\
	W^*
	}
$$  that is described by the $4 \times 
4$-pfaffians of the skew symmetric matrix
\[
	\begin{pmatrix}
	   0 & w^\perp_{12} & w^\perp_{13} & w^\perp_{14} & w^\perp_{15} \\
	   -w^\perp_{12} & 0 & w^\perp_{23} & w^\perp_{24} & w^\perp_{25} \\
	   -w^\perp_{13} & -w^\perp_{23} & 0 & w^\perp_{34} & w^\perp_{35} \\
	   -w^\perp_{14} & -w^\perp_{24} & -w^\perp_{34} & 0 & w^\perp_{45} \\
	   -w^\perp_{15} & -w^\perp_{25} & -w^\perp_{35} & -w^\perp_{45} & 0 \\
	\end{pmatrix} + t
	\begin{pmatrix}
	   0 & w^*_{12} & w^*_{13} & w^*_{14} & w^*_{15} \\
	   -w^*_{12} & 0 & w^*_{23} & w^*_{24} & w^*_{25} \\
	   -w^*_{13} & -w^*_{23} & 0 & w^*_{34} & w^*_{35} \\
	   -w^*_{14} & -w^*_{24} & -w^*_{34} & 0 & w^*_{45} \\
	   -w^*_{15} & -w^*_{25} & -w^*_{35} & -w^*_{45} & 0 \\
	\end{pmatrix}
\] and whose discriminant is isomorphic to the cone over the dual 
variety $\Edual$ of $E$. Homogeneous coordinates of the
$\PP^5$ are given by a basis of $\Wperp$ and $t$.
\end{cor}

\begin{proof} Let $w^*_{ij}$ be an nonzero element. Over the affine 
chart given by $w_{ij}^*=1$ the family
$\GGdualblowup \to \CC^4$ is described by
\[
	\begin{pmatrix}
	   0 & w^\perp_{12} & w^\perp_{13} & w^\perp_{14} & w^\perp_{15} \\
	   -w^\perp_{12} & 0 & w^\perp_{23} & w^\perp_{24} & w^\perp_{25} \\
	   -w^\perp_{13} & -w^\perp_{23} & 0 & w^\perp_{34} & w^\perp_{35} \\
	   -w^\perp_{14} & -w^\perp_{24} & -w^\perp_{34} & 0 & w^\perp_{45} \\
	   -w^\perp_{15} & -w^\perp_{25} & -w^\perp_{35} & -w^\perp_{45} & 0 \\
	\end{pmatrix} + t
	\begin{pmatrix}
	   0 & w^*_{12} & w^*_{13} & w^*_{14} & w^*_{15} \\
	   -w^*_{12} & 0 & w^*_{23} & w^*_{24} & w^*_{25} \\
	   -w^*_{13} & -w^*_{23} & 0 & w^*_{34} & w^*_{35} \\
	   -w^*_{14} & -w^*_{24} & -w^*_{34} & 0 & w^*_{45} \\
	   -w^*_{15} & -w^*_{25} & -w^*_{35} & -w^*_{45} & 0 \\
	\end{pmatrix}
\] where $t=0$ is the equation of the exceptional divisor that 
intersects all fibers in the elliptic curve $E_\perp$. Omitting
the condition $w_{ij}^*=1$ we obtain a family as in the proposition. 
By construction this is independent of our choice of
affine chart. Over $0 \in W^*$ its fiber is the projective closure of the cone over $E_\perp$. 
If $w \in W^*$ is a point different from $0$, $Y_w$ is the
fiber over this point, and $x=(w^\perp_{ij},t) \in Y_w$ is a singular 
point, then $x_\lambda = (w^\perp_{ij},\lambda^{-1}t)$
will be a singular point of $Y_{\lambda w}$ for $\lambda\not=0$. 
This proves that the discriminant of $X \to W^*$ is a cone.
With \ref{pDualDiscriminant} we obtain that it is 
isomorphic to the cone over the dual variety $\Edual$ of $E$.
\end{proof}


\begin{theorem} \label{theoA4} Let 
$\Gamma_\tau=\ZZ+\tau \ZZ$ be a lattice, $E=E_\tau=\CC/\Gamma_\tau$ the corresponding elliptic curve and $L(\tau)$ the matrix
$$L(\tau):= \left( \begin{array}{ccccc} 0 & x_6 & 0 & -x_4 & -2x_3\\[2ex] 
-x_6 & 0 & -\frac{2}{3}x_4 & -4x_3 & -8x_2\\[2ex] 0 &
\frac{2}{3}x_4 & 0 & -\frac{2}{3}x_2-\frac{1}{18}g_2(\tau)x_6 & 
-\frac{4}{3}x_1\\[2ex] x_4 & 4x_3 &
\frac{2}{3}x_2+\frac{1}{18}g_2(\tau)x_6 & 0 &
\frac{4}{3}g_2(\tau)x_4+2g_3(\tau)x_6\\[2ex] 2x_3 & 8x_2 & \frac{4}{3}x_1 & 
-\frac{4}{3}g_2(\tau)x_4-2g_3(\tau)x_6 & 0
\end{array} \right).$$ 

For a fixed $\tau$ the $4 \times 4$ Pfaffians of $L$ yield the cone $(X,0)$ over $E_\tau$ in its
fifth Weierstrass embedding as elliptic normal curve of degree $5$. If $N$ is the matrix
$$N:=\left( \begin{array}{ccccc} 0 & 0 & \y_5 & -3 \y_4 & 0\\[2ex] 0 
& 0 & 2 \y_4 & - 24 \y_3 & 0\\[2ex] - \y_5 & - 2 \y_4 & 0
& - 4 \y_2 & 0\\[2ex] 3 \y_4 & 24 \y_3 & 4 \y_2 & 0 & 48 \y_0\\[2ex] 
0 & 0 & 0 & -48 \y_0 & 0
\end{array} \right),$$
the family $L+N$ parametrized by $(\tau,\y_0, \y_2, \y_3, \y_4, \y_5)$ constitutes a versal deformation of $(X,0)$ and the matrix $A$ from \ref{exA} is a discriminant matrix for the discriminant of this deformation.

\end{theorem}

\begin{proof} 

Let $(y_0,y_2,y_3,y_4,y_5)$ be the coordinates of the vector space 
$W$, $(x_1,x_2,x_3,x_4,x_6)$ the coordinates of $\Wperp$,
and denote by
$(y^\ast_0, y^\ast_2, y^\ast_3, y^\ast_4, y^\ast_5)$ and $(x^\ast_1, 
x^\ast_2, x^\ast_3, x^\ast_4, x^\ast_6)$ the dual
coordinates of $W^\ast$ and $(\Wperp)^\ast$ respectively. The matrix
$$M:= \left( {\begin{array}{ccccc} 0 & g_3\,y_0 + \frac {1}{3}g_2y_2 
& y_5 &  - \frac {1}{3}g_2y_0  + \frac {2}{3}y_4  & y_3
\\[2ex] -g_3y_0 - \frac {1}{3}g_2y_2 & 0 & - \frac {1}{2}g_2y_0  - y_4 &  -
\frac{1}{2}y_3  &  -\frac{1}{2}y_2 \\[2ex]
  - y_5 & \frac {1}{2}g_2\,y_0  + y_4 & 0 & 6y_2 & 0 \\[2ex]
\frac {1}{3}g_2\,y_0  - \frac {2}{3}y_4  & \frac {1}{2}y_3  &  - 6y_2 
& 0 & - \frac {1}{2}y_0  \\[2ex]
  - y_3 & \frac {1}{2}y_2  & 0 & \frac {1}{2}y_0  & 0
\end{array}}
  \right) = (w_{ij})
$$ gives a mapping $\widetilde{M}: \bigwedge^2 V \to W$, $v_i \wedge 
v_j \mapsto w_{ij}$, and the matrix $L = (w_{ij}^\perp)$
a mapping $\widetilde{L}: 
\bigwedge^2 V^\ast \to \Wperp$, $v_i^\ast \wedge v_j^\ast
\mapsto w_{ij}^\perp$. The sequence
$$0 \longrightarrow W^\ast 
\stackrel{\widetilde{M}^T}{\longrightarrow} \bigwedge^2 V^\ast
\stackrel{\widetilde{L}}{\longrightarrow} \Wperp \longrightarrow 0$$ is exact. Indeed
$y_0^*$ is mapped under $\widetilde{M}^T$ to 
$$
\left( {\begin{array}{ccccc} 
0 & g_3 & 0 & -\frac{1}{3}g_2 & 0 \\
-g_3& 0 &-\frac{1}{2}g_2 & 0  & 0\\
0 & \frac{1}{2}g_2 & 0 & 0 & 0 \\
\frac{1}{3}g_2 & 0 & 0 & 0 & -\frac{1}{2} \\
0 & 0 & 0 & \frac{1}{2} & 0 
\end{array}} \right)
$$	
and evidently this is mapped to zero by $\widetilde{L}$. A similar calculation for the other basis elements shows that the sequence above is a complex. Furthermore the entries of $M$ and $L$ generate $W$ and $\Wperp$ respectively. Consequently both $\widetilde{M}$ and $\widetilde{L}$ have full rank $5$ and the sequence is exact.

We now claim that the elliptic curve $E=\PP(W) \cap \GG$  is 
parametrized by the Weierstrass embedding
$$\begin{array}{ccccc} \Phi & : & E & \to & \PP(W) \\ & & z & \mapsto 
& (1:\wp(z):\wp'(z):\wp''(z):\wp'''(z))
\end{array} .$$
Evaluating the $4 \times 4$-Pfaffians of $M$ on $\Phi(E)$, that is setting 
$y_0=1$ and $y_k = \wp^{(k-1)}$, $k=2,3,4,5$,
we must show:
\begin{eqnarray*}
\frac{1}{2} \wp'(z) \wp'''(z) -\frac{2}{3} \wp''(z)^2 + 2 g_2 
\wp(z)^2 + 6 g_3 \wp(z) + \frac{1}{6} g_2^2 & = & 0,\\
\wp'(z) \wp''(z) + \frac{1}{2} g_2  \wp'(z) - \frac{1}{2} \wp(z) 
\wp'''(z), & = & 0,\\
\frac{1}{2}\wp'(z)^2 - \frac{1}{3} \wp(z) \wp''(z) + \frac{1}{3} g_2 
\wp(z) + \frac{1}{2} g_3  & = & 0, \\ 6 \wp(z) \wp'(z) -
\frac{1}{2}  \wp'''(z) & = & 0, \\
\frac{1}{2}  \wp''(z) - 3 \wp(z)^2 + \frac{1}{4} g_2  & = & 0.
\end{eqnarray*}
These relations are a consequence of the following classical relations 
between the Weierstrass function and its
derivatives:
\begin{eqnarray*} (\wp')^2 & = & 4 \wp^3 - g_2 \wp -g_3, \\
\wp'' & = & 6 \wp^2 - \frac{1}{2} g_2, \\
\wp''' & = & 12 \wp \wp'.
\end{eqnarray*}
Note that \cite{Hu} presents the corresponding calculation 
for a different, but projectively equivalent embedding of an elliptic 
curve in $\PP^4$.

On the other hand, we consider the embedding $\Psi$ of the 
elliptic curve $E_\perp = \PP(\Wperp) \cap
\check{\GG}$ into $\PP(\Wperp)$ through
$$\Psi(z)= (- \frac{1}{24}  \wp'''(z):
\frac{1}{6}
 \wp''(z): - \frac{1}{2}  \wp'(z):
\wp(z):1) .$$ 
Setting
$$x_k = \frac{(-1)^{4-k}}{(5-k)!} \wp^{(4-k)}, \quad 
k=1,2,3,4,6;$$ 
with $\wp^{(-2)}\equiv 1$; this embedding is described by the $4 \times 4$ Pfaffians 
of the matrix $L$.

For a point
$(y^\ast_0, y^\ast_2, y^\ast_3, y^\ast_4, y^\ast_5) \in \PP(W^\ast)$
the corresponding $\PP^3 \subset \PP^4= \PP(W)$ is defined by
$$y^\ast_0y_0 + y^\ast_2y_2 + y^\ast_3y_3 + y^\ast_4 y_4 + y^\ast_5y_5 =0.$$
Substituting $y^\ast_k = \frac{(-1)^{k}}{(k-1)!} \y_k$, for $k=0,2,3,4,5$, 
the equation of $\PP^3 \cap E$ becomes the equation \ref{prop:Lambda}(\ref{eqLambda}).

The matrix
$$N=\left( \begin{array}{ccccc} 0 & 0 & \y_5 & -3 \y_4 & 0\\[2ex] 0 & 
0 & 2 \y_4 & - 24 \y_3 & 0\\[2ex] - \y_5 & - 2 \y_4 & 0
& - 4 \y_2 & 0\\[2ex] 3 \y_4 & 24 \y_3 & 4 \y_2 & 0 & 48 \y_0\\[2ex] 
0 & 0 & 0 & -48 \y_0 & 0
\end{array} \right) = (w_{ij}^\ast)$$ defines a  mapping 
$\widetilde{N} : \bigwedge^2 V^\ast \to W^\ast$, $v^\ast_i \wedge
v^\ast_j \mapsto w^\ast_{ij}$. One easily verifies that this is a 
left inverse of
$\widetilde{M}^T : W^\ast \to \bigwedge^2 V^\ast$ and  defines a 
splitting $\bigwedge^2 V^\ast \cong \Wperp \oplus W^*$.  

According to results of H.~Pinkham \cite{Pi1,Pi2} (see also \cite{M}), deformations of the cone $C(E_\perp)$ over the elliptic curve $E_\perp$ lift to projective deformations of the projective closure $\overline{C(E_\perp)}$ in $\PP^5$. Therefore
\ref{theoA4} follows from \ref{prop:discLambda} and 
\ref{CorHC}.
\end{proof}


\begin{remark}
The deformation of \ref{theoA4} is topologically trivial along the $\tau$-axis. This follows from the fact that the vector field $\partial / \partial \tau$ on $\HH \times \CC^5$ can be locally lifted to a vector field on $\CC^5 \times \HH \times \CC^5$. To see this, note that 
$$\frac{\partial g_2}{\partial \tau} = \frac{3}{\pi \sqrt{-1}} g_3, \quad \frac{\partial g_3}{\partial \tau} = \frac{1}{6 \pi \sqrt{-1}} g_2^2,$$
according to \cite{FS}, and that $g_2$ and $g_3$ do not vanish at the same time. From \cite{Ok} it follows that the complement of the discriminant in $\HH \times \CC^5$ is a $K(\pi,1)$-space.
\end{remark}



\end{document}